\documentclass[sts]{imsart}

\RequirePackage{amsthm,amsmath,amsfonts,amssymb,mathtools}
\RequirePackage[authoryear]{natbib}
\RequirePackage[colorlinks,citecolor=blue,urlcolor=blue]{hyperref}
\RequirePackage{graphicx}
\usepackage{txfonts}
 \usepackage{color}

 \newcommand{\commentout}[1]{}

\newcommand{\pspace}{\ensuremath{\Delta}}
\newcommand{\param}{\ensuremath{\delta}}

\date{\today}

\usepackage{bm,bbm}
\allowdisplaybreaks

\startlocaldefs
\newcommand{\seqe}[2]{S_{#1}^{#2}}
\newcommand{\prode}[2]{E_{#1}^{#2}}
\newcommand{\E}{\mathbb{E}}
\newcommand{\R}{\mathbb{R}}
\newcommand{\cB}{\mathbf{B}}
\newcommand{\cF}{\mathbf{F}}
\newcommand{\cG}{\mathbf{G}}
\newcommand{\cL}{\mathbf{L}}
\newcommand{\cP}{\mathbf{P}}
\newcommand{\cQ}{\mathbf{Q}}
\newcommand{\cV}{\mathbf{V}}
\newcommand{\cW}{\mathbf{W}}
\newcommand{\Tau}{\mathcal{T}}

\endlocaldefs

\begin{document}

\begin{frontmatter}

\title{Game-Theoretic Statistics and \\
Safe Anytime-Valid Inference}
\runtitle{Game-theoretic statistics and safe anytime-valid inference}


\begin{aug}
\author[A]{\fnms{Aaditya} \snm{Ramdas}\ead[label=e1]{aramdas@cmu.edu}},
 \author[B]{\fnms{Peter} \snm{Gr\"unwald}\ead[label=e2]{pdg@cwi.nl}},
 \author[C]{\fnms{Vladimir} \snm{Vovk}\ead[label=e3]{v.vovk@rhul.ac.uk}}
 \and
 \author[D]{\fnms{Glenn} \snm{Shafer}\ead[label=e4]{gshafer@business.rutgers.edu}}
 \address[A]{Carnegie Mellon University, Pittsburgh, USA \printead{e1}.}
 \address[B]{Centrum Wiskunde \& Informatica, Amsterdam, Netherlands \printead{e2}.}
 \address[C]{Royal Holloway London, UK \printead{e3}.}
 \address[D]{Rutgers University, USA \printead{e4}.}
\end{aug}

\begin{abstract}
Safe anytime-valid inference (SAVI) provides measures of statistical evidence and certainty---e-processes for testing and confidence sequences for estimation---that remain valid at all stopping times, accommodating continuous monitoring and analysis of accumulating data and optional stopping or continuation for any reason. These measures crucially rely on test martingales, which are nonnegative martingales starting at one. Since a test martingale is the wealth process of a player in a betting game, SAVI centrally employs game-theoretic intuition, language and mathematics. We summarize the SAVI goals and philosophy, and report recent advances in testing composite hypotheses and estimating functionals in nonparametric settings.
\end{abstract}

\tableofcontents


\end{frontmatter}

\section{Introduction}

Stop when you are ahead. Increase your bet to make up ground when you are behind.  This is called martingaling in the casino.  It often succeeds in the short or medium term, leading novice gamblers to think they can beat the odds and day traders to think they can beat the market \citep{Dimitrov/etal:2022}.  The same delusion arises in science, where sampling until a significant result is obtained is an important source of irreproducibility.

The fallacy of sampling until a significant result is obtained has been discussed by statisticians at least since the 1940s, when \cite{Feller:1940} saw it happening in the study of extra-sensory perception. \cite{Anscombe:1954}  famously called it ``sampling to a foregone conclusion'', and this inevitability was also pointed out by~\cite{robbins1952some}. 

But disapproval by statisticians has hardly dented the prevalence of the practice.  In one widely publicized example, a team of researchers apparently demonstrated benefits from ``power posing'' \citep{Carney/etal:2010}.  The lead author later disavowed the conclusion and identified the team's peeking at the data as one of her reasons \citep[Fact 5]{Carney:2022}: 
\begin{quote} 
We ran subjects in chunks and checked the effect along the way. It was something like 25 subjects run, then 10, then 7, then 5. Back then this did not seem like p-hacking. It seemed like saving money (assuming your effect size was big enough and p-value was the only issue).  
\end{quote}
Ten years ago, an anonymous survey of over 2000 psychologists found $56\%$ admitting to ``deciding whether to collect more data after looking to see whether the results were significant'' \citep{john2012measuring}.

Bayesian inference with a prior defined by a statistician's beliefs before seeing any of the data is not affected by (planned) peeking.  Problems quickly arise, however, when default or pragmatic priors are used to test composite null hypotheses.  These problems are especially severe for commonly used pragmatic priors that depend on the sample size, covariates, or other aspects of the data \citep{deHeide/Grunwald:2021}.  

As emphasized by~\cite{johari2022always,howard2021time,GrunwaldHK19,Shafer:2021,pace2020likelihood}, amongst others,
we need to go beyond disapproval of peeking, and we instead should give researchers tools to fully accommodate it. The branch of mathematical statistics that enables this, sequential analysis, was brilliantly launched in the 1940s and 1950s by  Wald, Anscombe,  Robbins, and others. The innovations introduced by Robbins, Darling, Siegmund and Lai included {\em confidence sequences\/} that are valid at any and all times and {\em tests of power one}. But these ideas occupied only a small niche in sequential analysis research until around 2017. 
Since then, interest has exploded and much conceptual progress has been made in parallel threads, which we attempt to summarize. 

This new methodology differs from traditional statistical testing in the way it quantifies evidence against statistical hypotheses.  The traditional approach casts doubt on a hypothesis when a selected test statistic takes too extreme a value.  This leads to quantifying evidence against the hypothesis by the \textit{p-value}---the probability the hypothesis assigns to the test statistic being so large.  The new methodology instead casts doubt on a hypothesis when a selected nonnegative statistic is large relative to its expected value.  Imagining that we bought the statistic for its expected value when we selected it, we call the ratio of its realized to its expected value a \textit{betting score} and take this as a measure of our evidence.  In the case of a composite hypotheses, we use the infimum of betting scores for the multiple hypotheses and call this an \textit{e-value}.  The sequential analog is an \textit{e-process}---a sequence of e-values that monitor the accumulation of evidence.  E-processes permit anytime-valid inference; we can repeatedly decide whether to collect more data based on the current e-value without invalidating later assessments, stopping whenever and for any reason whatsoever.  This anytime-validity is a form of \textit{safety}. This safety may come with a price, of course;  there may (or may not) be tradeoffs between safety and power;  see Section~\ref{app:price}.

From a technical point of view, the new methodology is based on the concept of a test martingale, along with its betting interpretation.  Although martingales became important in probability theory more than a half-century ago, their potential has still not been fully exploited in statistics, and the new emphasis on \emph{nonnegative} (super)martingales has produced a plethora of powerful new methods.  These include confidence sequences for many functionals that can be used with multi-armed bandits and new sequential tests for composite null hypotheses.  This responds to the need for rigorous methods in settings that have emerged with the development of information technology in the past half-century, including ``living meta-analysis''~\citep{Schure/Grunwald/Ly:2021}, the industrial use of A/B testing~\citep{johari2022always} and bandit experiments~\citep{howard2022sequential}.

The new methods can be most clearly presented in the language of game-theoretic probability \citep{Shafer/Vovk:2001,Shafer/Vovk:2019}. Here successive observations are Reality's moves in a game.  Two other players move before Reality on each round:  Forecaster gives probabilities for the outcome, and Skeptic bets by choosing a real-valued function of the outcome, paying its expected value, and receiving its realized value.  If Skeptic always chooses nonnegative functions, then the factor by which he multiplies his money (the ratio of the realized to the expected value under forecaster's probabilities) is his ``betting score'' or ``e-value'' \citep{Shafer:2021}.  If he reinvests his money on each round, the betting scores multiply, producing cumulative betting scores that are products of the betting scores for each round so far.  Because Skeptic is a free agent, the option of stopping or continuing or even switching to a different experiment on the next round is intrinsic to the game, and the cumulative betting score or e-value quantifies the evidence against the Forecaster (and his probabilities): Skeptic refutes the odds by making money betting at those odds; more money is more evidence that the odds do not reflect reality.

Betting games often fit statistical practice better than measure-theoretic probability models.  In particular, they accommodate fully the opportunistic behavior that we want to allow.  George Barnard, in his review of Wald's book on sequential analysis \citep{Barnard:1947}, called for embedding statisticians in the sequential decision-making of experimental scientists, in which each batch of observations is followed by deliberation about whether to stop or to continue, perhaps with a modified experiment.  The use of a prespecified stopping time, which prescribes continuing only until a certain data-dependent condition is met, obscures or erases this sequential deliberation, pretending that all the decisions flow from a stopping strategy adopted in advance. %
 Barnard's suggestion is better captured by our game-theoretic framework, where a single stopping rule is replaced by notions of evidence that remain valid at \emph{any} stopping time not specified in advance.

Because most readers will be unfamiliar with game-theoretic probability as developed by \citet{Shafer/Vovk:2001,Shafer/Vovk:2019}, we use the relatively familiar apparatus of measure theory (filtrations, stopping times, martingales, etc.) and new concepts defined within that apparatus (e-values, e-processes, etc.). Frequently, however, we return to the betting story, where our martingales are wealth processes for Skeptic. 

\subsection{Basic terminology}

We begin with a sample space $\Omega$ equipped with a filtration $\cF \equiv ({\cF_t})_{t\geq 0}$ (an increasing nested sequence of $\sigma$-fields), and a set $\Pi$ of probability distributions on $(\Omega, \cF)$. We assume that some distribution $P \in \Pi$ governs our data $X \equiv (X_1,X_2,\dots)$. The variables $X_1,X_2,\dots$ need not be independent and identically distributed (iid) under $P$. We use $X^t$ as a shorthand for $X_1,\dots,X_t$.

When we say we are testing $\cP$, we mean that we are testing the null hypothesis $H_0$ that 
$P \in \cP$.  When we say we are testing $\cP$ against $\cQ$, we mean that the alternative hypothesis $H_1$ is that $P \in \cQ$. Typically, $\cP$ and $\cQ$ are either non-intersecting or nested subsets of $\Pi$. We always use boldface $\cP,\cQ$ for sets of distributions, and a normal $P,Q$ for a single distribution.

A sequence of random variables $Y \equiv (Y_t)_{t\geq 0}$ is called a \textit{process} if it is adapted to $\cF$---i.e., if $Y_t$ is measurable with respect to $\cF_t$ for every $t$. 
Often $\cF_t := \sigma(X^t)$, with $\cF_0$ being trivial ($\cF=\emptyset,\Omega\}$), and in this case $Y_t$ being measurable with respect to $\cF_t$ means that $Y_t$ is a measurable function of $X_1,\dots,X_t$.  
But $\cF$ is sometimes a coarser filtration (we discard information, see e.g.\ Section~\ref{sec:ttest}) or a richer one (we add external randomization).%
\footnote{The filtration may be coarsened, as explained by Alan Turing (\citeyear[p.1]{Turing:1941}):  ``When the whole evidence\dots
is taken into account it may be extremely difficult to estimate the probability of the event, \dots
may be better to form an estimate based on a part of the evidence \dots''
}
$Y$ is called \emph{predictable} if $Y_t$ is measurable with respect to $\cF_{t-1}$.

A stopping time (or rule) $\tau$ is a nonnegative integer valued random variable such that $\{\tau \leq t\} \in \cF_t$ for each $t \geq 0$. In words: we know at each time whether the rule is telling us to stop or keep going. Denote by $\Tau$ the set of all stopping times, including ones that may never stop.

\subsection{A terse technical summary of the paper} 

We give a short technical summary below, foreshadowing topics to be defined and discussed in more depth later. 

The field of safe anytime-valid inference (SAVI) aims to develop measures of statistical evidence and certainty that remain valid at arbitrary stopping times (possibly unknown in advance), accommodating continuous monitoring and analysis of accumulating data and optional stopping or continuation for any reason. There is a strong sense in which \emph{admissible} SAVI methods --- power-one tests, confidence sequences, anytime-valid p-values and e-processes --- \emph{must} rely centrally on nonnegative martingales~\citep{ramdas2020admissible}. Nonnegative (super)martingales are endowed with a strong and direct connection to gambling: every nonnegative supermartingale corresponds to a wealth process in some game, and vice versa (every ``fair/legal'' gambling strategy to test the null hypothesis results in a wealth process that is a nonnegative supermartingale).

These facts give rise to the following central principle in game-theoretic statistics: ``testing by betting''.  In order to test a null hypothesis $\cP$ against an alternative $\cQ$, we set up a game such that (a) if the null is true, meaning $P \in \cP$, then no betting strategy can reliably make money (any gambler's wealth is a nonnegative supermartingale), and (b) if the null is false, meaning $P\in \cQ$, it is possible to bet smartly to make money in that game. This principle arguably has roots dating back (at least) to~\cite{Ville:1939}, and was recently discussed in depth  in the point null case by~\cite{shafer2011test} and~\cite{Shafer:2021} and for composite nulls by~\cite{GrunwaldHK19,waudby2020estimating}, etc.

The game works as follows.
Before observing $X_t \in \mathcal X$,  Skeptic puts forward a bet $S_t:\mathcal X \to [0,\infty]$, which satisfies 
\begin{align}\label{eq:cond-bet}
\E_P[S_t(X_t)] \leq 1 \text{ for every } P\in\cP.
\end{align}
Then, $X_t$ is revealed. 
%
%
The interpretation is that at each time $t$, one can buy, at the price of 1 monetary unit, a ticket that will pay off $S_t(X_t)$ units. One can buy as many tickets as one likes. (\ref{eq:cond-bet}) simply expresses that under the null, one does not expect to get back more than one invests in this game. At time 1, Skeptic invests 1 monetary unit; at each time $t$, she reinvests all the money she observed so far. Skeptic's wealth after $t$ steps is then clearly   given by $\prod_{i=1}^t S_i(X_i)$, which is nonnegative by definition, and 
easily checked to be a supermartingale under $P$. 
%
%
If $\cQ$ is appropriately separated from $\cP$,  good betting strategies can force the wealth in setting (b) (alternative is true) to grow to infinity exponentially fast, and we wish to maximize the exponent. Maximizing the exponent corresponds to maximizing the expected logarithm of the wealth; such a ``log-optimality'' objective has information-theoretic roots dating back to \cite{Kelly56} and \cite{breiman1961optimal}, but also (implicitly or explicitly) appears in the work of Ville, Wald, Robbins, etc. 

For testing a point null $P$ against a point alternative $Q$, the log-optimal bet is simply given by the likelihood ratio $S_t = dQ_t/dP_t$, where $P_t$ and $Q_t$ are the conditional probabilities for the $t$th observation given the past, under $P$ and $Q$ respectively. Thus the realized likelihood ratio of $Q$ against $P$ is precisely the optimal wealth of a gambler betting against $P$. This central fact provides much intuition for extensions and generalizations.

For composite alternatives $\cQ$, the Skeptic often hedges their bets by not betting all their money on a single $Q \in \cQ$, instead spreading their investment over $\cQ$ using a mixture (``prior'') distribution $R$.
To illustrate using the toy case in which $\cQ$ is countable and $R$ has mass function $r$, $r(q) = a$ would mean that a fraction $a$ of Skeptic's money is invested in $q$ --- importantly in general $R$ neither has a frequentist (`drawing from an urn') nor a Bayesian (`belief')  interpretation here.
This \emph{method of mixtures} plays a central role in this paper: an instance of Laplace's method for approximating a maximum by an integral, 
it appears directly within our anytime-valid context in~\cite{Ville:1939} and~\cite{robbins1970statistical}, and in broader sequential contexts
in \cite{Wald:1947} and~\cite{cover1974universal}, among many others.  

The most interesting questions in this area involve composite (and often nonparametrically specified) nulls $\cP$.
Indeed, there really was no general theory for dealing with composite nulls until 2017 ---  when, almost out of the blue, several generic proposals for dealing with composite nulls appeared. Arguably, it is this development which caused the aforementioned explosion of interest in the area --- suddenly there was an indication that eventually almost any interesting statistical testing or estimation problem could be converted into an anytime-valid version with a gambling interpretation. 

For such composite $\cP$, a fascinating phenomenon sometimes presents itself: for some ``extremely rich'' nulls $\cP$, the game described above is hopelessly restraining: the constraint~\eqref{eq:cond-bet} is too stringent, and the only functions $S_t$ that satisfy it are either constant or decreasing (meaning that they cannot increase under any alternative). This happens, for example, when testing exchangeability or testing log-concavity; see Section~\ref{sec:exch+LC} for references and details.

Luckily, generalizing the above game protocol resuscitates the approach.
There appear to be two different types of generalized games: (a) one can restrict the amount of information available to the Skeptic by introducing a third player (an ``Intermediary'') who throws away some information revealed by Reality (mathematically, Skeptic operates in a shrunk filtration), (b) one can instead make the Skeptic play many games in parallel, each against a different subset of $\cP$, with the Skeptic's net wealth being their worst wealth across all the parallel games. In the first case, Skeptic's wealth may remain a nonnegative (super)martingale, but in the second case, their wealth is an \emph{e-process} (under the null, their wealth is upper-bounded by a different supermartingale in each game, and thus is bounded by one at any stopping time). While these solutions may seem almost magical at first glance, they both yield fruit for the same problem mentioned above of testing exchangeability: approach (a) is used in \citet[Part~III]{Vovk/etal:2022} and approach (b) in \citet{ramdas2022testing}. The latter work, along with~\cite{ruf2022composite}, together show the centrality of e-processes in game-theoretic statistics: e-processes exist for many $\cP$ for which nonnegative (super)martingales do not.

When $\cP$ and $\cQ$ have a common reference measure, meaning that likelihood-ratio based methods are still in play, two key ideas stand out:  universal inference~\citep{wasserman2020universal}, and the reverse information projection~\citep{GrunwaldHK19}. 
The former always yields an e-process, but latter always results in an e-value which can be multiplied across batches of data to yield a supermartingale. But \emph{sometimes} the latter also directly yields an e-process (and when it does, it dominates universal inference).

When $\cP$ and $\cQ$ do not have any common reference measures --- and thus likelihood-ratio based methods may not make any sense at the outset --- the \emph{design} of nonnegative (super)martingales or e-processes occupies center stage. Sometimes, the nonparametric definition of $\cP$ directly yields a natural game, like when testing if a ``subGaussian'' mean is positive~\citep{Darling/Robbins:1967b}. Other times, one must design new games in possibly shrunk filtrations, which may not be obvious at the outset, like in two-sample testing~\citep{shekhar2021game}. 

The entire discussion above was centered on testing by betting, because this typically forms the technical heart of other problems that are not cast explicitly as testing. For example, appropriate duality concepts and inversions allow us to translate many of these results into those for estimation of appropriate functionals using confidence sequences (Section~\ref{sec:nonparametric}). Both e-processes and confidence sequences can in turn be extended to other problems like change detection (Section~\ref{sec:change}), model selection, etc. 

In fact, our investigations reveal a curious phenomenon: at the \emph{heart} of many (and plausibly, \emph{all}) nonparametric testing and estimation problems 
%
is a ``hidden'' game (often not unique: the same $\cP$ and $\cQ$ may be associated with different filtrations and betting strategies that are e-processses under $\cP$ and make money under $\cQ$).  Further, explicitly bringing out 
such games (and betting well in them) can yield powerful new methodology as well as new theoretical insights~\citep{howard2020time}

A full understanding of when and why this happens is open, but we provide one hint here. Likelihood ratios have been at the center of statistics for nearly a century. Nonnegative (super)martingales and e-processes are simply nonparametric, composite generalizations of likelihood ratios, and these have been found to exist in dozens of problems where one cannot even begin to talk about likelihood-ratio based methods. Thus, these tools give us a way to work implicitly with likelihood ratios, even when there appears to be no explicit way to do so. Given the power (and sometimes optimality) of likelihood-ratio based tests in parametric settings, we perhaps get a hint of the power of our game-theoretic approaches in composite (often nonparametric) settings.

The rest of this paper will formally define the key concepts, and provide  technical details of the aforementioned methods and phenomena in different problem settings.

\section{Central concepts}

In the sequel, we leave measurability assumptions and other measure-theoretic details implicit so far as possible.

\subsection{E-values}\label{sec:evalues}

An \emph{e-variable} for $\cP$ is a nonnegative random variable $E$ such that $\E_P[E]\leq 1$ for all $P \in \cP$. Its realized value, after observing the data, is an \textit{e-value}.\footnote{Observe that we use boldface $\E$ for expectation and normal $E$ for e-values. The ``e'' in e-value stands both for ``evidence'' (because it quantifies statistical evidence against the null) and  for ``expectation'' (because its central property is its expectation).} Often we call $E$ itself an e-value, blurring the distinction between the random variable and its realized value.  (The term ``p-value'' is also often used for both random variables and their values.)

When $\E_P[E] = 1$, we call the e-value $E$ a \textit{unit bet against $P$}. This name evokes a story in which expected values are prices of payoff: the Forecaster predicts that $X\sim P$, and in order to bet against them, a Skeptic could buy one unit of $E$, for the price of $1$, delivering the Skeptic a payoff of $E(X)$. Mathematically, a unit bet against $P$ is simply\footnote{In some sense, statisticians have always been using e-values (and test martingales), because likelihood ratios are the most important example of e-values (and test martingales). But this direct analog only holds when testing a single distribution $P$. The power and utility of e-values, test (super)martingales and e-processes are truly realized only dealing with a composite (and sometimes nonparametric) $\cP$.
} a likelihood ratio $dQ/dP$ for some alternative $Q$.  This is elementary when we use probability densities:
\begin{itemize}
\item $\E_P[E] = 1$ can be written as $\int E(x)p(x) dx$=1, so that $q:=E\times p$ is a density, and $E=q/p$.
    \item If $q$ and $p$ are $Q$'s and $P$'s densities, then  $\E_P[q/p]=\int p(x)\frac{q(x)}{p(x)}d x = 1$.
\end{itemize}

We use e-values when data are treated as a batch.  Their dynamic counterparts are test martingales and e-processes, introduced next.

\subsection{Test (Super)Martingales}\label{subsec:test}

A process $M$ is a \emph{martingale} for $P$ if 
\begin{equation}\label{eq:martingale}
    \E_P[M_t \mid \cF_{t-1}] = M_{t-1} 
\end{equation}
for all $t\geq 1$. $M$ is a \emph{supermartingale} for $P$ if it satisfies \eqref{eq:martingale} with ``$=$'' relaxed to ``$\leq$''.  A (super)martingale is called a \emph{test (super)martingale} if it is nonnegative and $M_0=1$.

Game-theoretically, a test martingale for $P$ is the wealth process of a gambler who bets against $P$. If $M$ is a test martingale for $P$, then $\E_P[M_t] = 1$ for any $t \geq 0$, and thus each $M_t$ is itself a unit bet against $P$; it is the factor by which $M$ multiplies its money from time $0$ to time $t$. Similarly, the optional stopping theorem implies that for \emph{any} stopping time $\tau$ --- even potentially infinite ---  $\E_P[M_\tau] \leq 1$, and thus each $M_\tau$ is also an e-value for $P$. 


The correspondence between unit bets against $P$ and likelihood ratios with $P$ as the denominator extends to a related correspondence for test martingales for $P$.  If $Q$ is absolutely continuous with respect to $P$, we can write
\begin{equation}\label{eq:ele}
   \frac{q(X^t)}{p(X^t)} =  \frac{q(X_1)}{p(X_1)} \frac{q(X_2\mid X_1)}{p(X_2\mid X_1)} \cdots \frac{q(X_t\mid X^{t-1})}{p(X_t\mid X^{t-1})},
\end{equation}
where $X^t:=(X_1,\dots,X_t)$, $p(X^t)$ is $P$'s density for $X^t$, and $q(X^t)$ is $Q$'s density for $X^t$.  Denote the sequence defined by \eqref{eq:ele} as $M$; then $M$ is a test martingale for $P$, and 
\begin{align}
\label{eq:tmart}
   M_t &= \prod_{i=1}^t B_i 
       = \frac{q(X^t)}{p(X^t)},\\
\label{eq:unit}
     \text{where ~}~ B_t &:= \frac{q(X_t\mid X^{t-1})}{p(X_t\mid X^{t-1})}.
\end{align}
Note that each $B_t$ is a unit bet against $P$, conditional on $\cF_{t-1}$; we call $B_t$ \textit{$M$'s unit bet on round $t$}.

Test martingales for $P$ are always of the form \eqref{eq:tmart}.  So choosing a test martingale for $P$ comes down to choosing an alternative $Q$.  In applications, constructing a test martingale for $P$ usually amounts to constructing the numerator $q(X_t\mid X^{t-1})$ in \eqref{eq:unit}; see Section~\ref{sec:compositeH1}. Test supermartingales can also be decomposed in the style of \eqref{eq:tmart}, where the $B_t$ are {\em single-round\/} e-values (i.e.\ defined as function on a single outcome $X_t$) conditional on $\cF_{t-1}$.

Test martingales 
become more interesting objects in the composite setting.

\subsection{Composite Test (Super)Martingales}

A process $M$ is a test (super)martingale for $\cP$ if it is a test (super)martingale for every $P \in \cP$. Such composite test (super)martingales are important in this paper. Composite test martingales also decompose as in~\eqref{eq:tmart}: for every $P \in \cP$, there is a $Q$ that is absolutely continuous with respect to $P$ and satisfies $M_t = q(X^t)/p(X^t)$; see~\citet[Proposition 4]{ramdas2020admissible}. In other words, \emph{composite test martingales are simultaneous likelihood ratios}.

Trivially, the constant process $M_t=1$ is a test martingale for any $\cP$, and a decreasing process is a test supermartingale for any $\cP$. We call a test (super)martingale \emph{nontrivial} if it is not always a constant (or decreasing) process. In particular, we would like test martingales for $\cP$ that increase to infinity under the alternative $\cQ$. But there may be no nontrivial test martingales if $\cP$ is too large.  In this case there may still be nontrivial test supermartingales (Section~\ref{sec:subgaussian}), but even these may not exist (Section~\ref{sec:exch+LC}). For this reason, we also need e-processes.

\subsection{E-processes}

A family $(M^P)_{P \in \cP}$ is a \emph{test martingale family} if $M^P$ is always a test martingale for $P$. 
A nonnegative process $\prode{}{}$ is called an \emph{e-process} for $\cP$ if there is a test martingale family $(M^P)_{P \in \cP}$ such that 
\begin{equation}
\label{def:eprocess1}
\prode{t}{} \leq M^P_t \text{ for every } P \in \cP,  t \geq 0.
\end{equation}
This type of definition was used by~\cite{howard2020time}, who used the name ``sub-$\psi$ process''.
In parallel, \cite{GrunwaldHK19} implicitly defined an e-process for $\cP$, also without using the name ``e-process'',  as a nonnegative process $\prode{}{}$ such that
\begin{equation*}
\E[\prode{\tau}{}] \leq 1 \text{ for every } \tau \in \Tau,P \in \cP.
\end{equation*}
In words, $\prode{}{}$ must be an e-value at any stopping time.
\cite{ramdas2020admissible} proved that the two definitions are equivalent and that if $\cP$ is ``locally dominated'', then \emph{admissible}\footnote{An e-process $E \equiv (E_t)_{t \geq 1}$ for $\cP$ is \emph{inadmissible} if there exists another e-process $E'$ for $\cP$ such that $E' \geq E$ ($E'_t \geq E_t$ almost surely $P$, for all $P\in\cP$ and all $t \geq 1$), and $E'_t > E_t$ with positive probability under some $P \in \cP$ and some $t \geq 1$; $E$ is admissible if it is not inadmissible.} e-processes (see Section~\ref{subsec:admissibility}) must satisfy 
\begin{equation}\label{eq:einf}
  \prode{t}{} = \inf_{P \in \cP} M^P_t
\end{equation}
for some test martingale family $(M^P)_{P \in \cP}$.
(Technically, the $\inf$ above is an ``essential infimum''.)

Whereas a test martingale for $P$ is the wealth process of a gambler who bets against $P$, an e-process for $\cP$ reports the minimum wealth across many simultaneous betting games, one against each $P \in \cP$, all with the same outcomes $X_1,X_2,\dots$ \citep[Section 5.4]{ramdas2022testing}.

The evidence against a null hypothesis as measured by a test martingale or e-process may decrease as we collect more data (indeed, gamblers may lose money as they  play a game, even if the odds are stacked in their favor). In order to obtain a measure of evidence that does not decrease with time, one can calculate the running maximum of the e-process $(\sup_{s \leq t} E_s)$ --- which is not an e-process --- and then adjust it back to being an e-process using a lookback calibrator; see~\citet{shafer2011test,dawid2011insuring} and~\citet[Section 4.7]{ramdas2022testing}.

\subsection{Ville's Theorem and Ville's Inequality}

The notion of a test martingale was first formulated by~\cite{Ville:1939}, though he simply called it a martingale.

Ville gave a proof, valid for any discrete-time stochastic process $P$, that an event $A$ has measure zero under $P$ if and only if there is a betting strategy that bets against $P$ and becomes infinitely wealthy if $A$ happens---i.e., a test martingale for $P$ that grows to infinity on all of $A$.
Moreover, $P(A) < \epsilon$ if and only if there is a test martingale for $P$ that exceeds $1/\epsilon$ on all of $A$.%
\footnote{\cite{Shafer/Vovk:2019} turn this around into a \emph{definition} of probability in the betting game.}
These results have been called \textit{Ville's theorem}  \cite[Section 9.1]{Shafer/Vovk:2019}. \cite{ruf2022composite} generalize Ville's theorem to composite $\cP$, but this cannot be accomplished by a test martingale, requiring e-processes instead.

Ville also showed that if $M$ is a test martingale for $P$, then for any $\alpha \ge 1$,
\begin{equation}\label{eq:ville1}
     P\left(\sup_t M_t \ge \alpha \right) \le \frac{1}{\alpha}.
\end{equation}
Ville called this the theorem of gamblers' ruin; a gambler who begins with unit capital and keeps betting until he wins the casino's entire capital $\alpha$ has little chance of succeeding.  More recently, the theorem has become known as \textit{Ville's inequality}.  For a self-contained proof see \cite{howard2020time} or \cite{Crane/Shafer:2020}.

Ville's inequality holds with equality for continuous-path (hence continuous-time) test martingales with unbounded total variation: a classic example would be the process $M_t = \exp(\lambda B_t - t\lambda^2/2)$, where $B_t$ is a standard Brownian motion and $\lambda$ is any nonzero constant. In discrete time, the main source of looseness is from ``overshoot'': continuous path martingales equal $1/\alpha$ at the instant of crossing it (they do not overshoot), but discrete time (and thus discrete path) test martingales are typically strictly larger than $1/\alpha$ at the first time of crossing; this causes some looseness. In practice (which is always in discrete time), Ville's inequality is ``quite tight'' and overshoot is often considered a second order effect.

Ville's inequality extends to statements about composite $\cP$: if $\prode{}{}$ is an e-process for $\cP$, then for every $\alpha \in (0,1)$,
\begin{equation}\label{eq:ville2}
\sup_{P \in \cP} P(\exists t \geq 1: \prode{t}{} \geq 1/\alpha) \leq \alpha.
\end{equation}
Equivalently, by \citet[Lemma 3]{howard2021time}, 
\begin{equation}\label{eq:ville3}
P(\prode{\tau}{} \geq 1/\alpha) \leq \alpha \text{ for every $\tau\in\Tau,P\in\cP$}.
\end{equation}
Ville's inequality plays a central role in converting e-processes into sequential tests or confidence sequences.  

Recently, \cite{wang2023extended} extended Ville's inequality to apply to \emph{nonintegrable} nonnegative supermartingales, such as those obtained when mixing likelihood ratios with an improper prior.

\subsection{Sequential Tests and their Families}
\label{sec:seq}
We consider test martingales and e-processes bona fide measures of evidence, with no need for thresholding. But we may want to make a binary decision based on this evidence. 
We define a (one-sided) sequential test in terms of rejection decisions like $(0,0,0,0,1,1,1,1,\dots)$, where a $0$ means that there is not yet enough evidence to reject the null, and a $1$ means that there is. In this formalization, a level-$\alpha$ sequential test $\psi \equiv (\psi_t)_{t \geq 1}$ is an increasing process consisting of 0-1 random variables 
such that 
\begin{equation}
    \label{eq:sequential}
P(\exists t\geq 0: \psi_t = 1) \leq \alpha \text{ for all $P \in \cP$. }
\end{equation}
\citet[Lemma 3]{howard2021time} proved that an \emph{equivalent} definition, with optional stopping made more explicit, is 
\[
P(\psi_\tau =1) \leq \alpha \text{ for any $\tau \in \Tau, P \in \cP$}. 
\]
%

It is easy to obtain a sequential test from a test martingale or e-process: simply reject the null (and stop) the first time the process reaches or exceeds $1/\alpha$. Indeed, Ville's inequality implies that
$\psi_t := \mathbf{1}(\sup_{s \leq t} M_s \geq 1/\alpha)$ is a sequential test.
We call a family $(\psi^P)_{P \in \cP}$, where $\psi^P$ is a sequential test for $P$, a \textit{sequential test family}.

Note that our sequential tests are different from Wald's original sequential tests, which dominated the area for a long time --- in these, the null hypothesis may finally be accepted and \emph{the stopping rule is specified beforehand}, determined by the desired Type-I/II error bounds $\alpha$ and $\beta$. Our framework resembles the ``power-one tests'' of~\cite{darling1968some}, where we do not commit to a stopping rule, and could keep going if we do not reject $\cP$.

\subsection{Anytime-valid p-values}

A random variable $p$ is a p-value for $\cP$ if $P(p \leq u) \leq u$ for all $P\in\cP$ and $u\in[0,1]$.  Like the e-value, this is a static concept.  
Anytime-valid p-values are the dynamic counterparts of p-values.

An \textit{anytime-valid p-value}~\citep{johari2022always,howard2021time} for $\cP$ is a process $p := (p_t)_{t \geq 1}$ such that $P(p_\tau \leq u) \leq u$ for any $\tau \in \Tau, P \in \cP,u\in[0,1]$. Equivalently, $P(\inf_t p_t \leq u) = P(\exists t\geq 1: p_t \leq u) \leq u$. In other words, with probability at least $1-u$ an anytime-valid p-value will never drop below $u$. So decisions to stop an experiment or to continue to collect data based on the current value of an anytime-valid p-value are \textit{safe}; they will not violate type-I error control. 

It is easy to check that if $M$ is an e-process for $\cP$ then $1/(\max_{s \leq t} M_s)$ is an anytime-valid p-value for $\cP$. In our framework, test martingales and e-processes are central objects for testing, and sequential tests and anytime-valid p-values take a secondary and derivative role.

\subsection{Confidence Sequences}
\label{sec:confidence}
When estimating some property of a distribution, like a mean (or a median), we  think of it as a functional $\phi: \Pi \to \Theta$ for some space $\Theta$, which is often a subset of $\mathbb{R}^d$. 

A $(1-\alpha)$-\textit{confidence sequence} (CS) is a sequence $(C_t)_{t \geq 0}$ of sets $C_t \subseteq \Theta$ such that 
\[
P(\forall t \geq 1: \phi(P) \in C_t) \geq 1-\alpha \text{ for all $P \in \Pi$}.
\]
As before, \citet[Lemma 3]{howard2021time} implies that a mathematically equivalent definition is to require
\[
P(\phi(P) \in C_\tau) \geq 1-\alpha \text{ for all } \tau \in \Tau, P \in \Pi.
\]
This dynamic concept can be contrasted with the concept of a confidence set (or interval).  A $(1-\alpha)$ confidence set, as usually defined, is required only to contain $\phi(P)$ with probability $1-\alpha$ for a sample of a fixed size or at a single fixed stopping time rather than at all stopping times. Confidence sequences remain valid under continuous monitoring (or peeking) and optional stopping, but confidence sets require the sample size or the stopping time to be fixed in advance of seeing any data.

One can construct a confidence sequence by inverting a family of sequential tests, or thresholding a test martingale family $(M^P)_{P \in \Pi}$:
\(
C_t := \{\phi(P): P \in \Pi, M^P_t <1/\alpha \}.
\)
Sometimes it is easier to construct a test martingale family $(M^\theta)_{\theta \in \Theta}$, where $M^\theta$ is a test martingale for $\{P: \phi(P) = \theta\}$. In that case, we would define
\begin{equation}
\label{eq:general-CS-by-inversion-2}
C_t := \{\theta \in \Theta:  M^\theta_t <1/\alpha \}.
\end{equation}

\subsection{Averaging e-values}\label{subsec:avg-e}
We can average e-values. By the linearity of expectations, if $E_1$ and $E_2$ are e-values for $\cP$, then $(E_1 + E_2)/2$ is as well, even if $E_1$ and $E_2$ are dependent (for example, calculated in different ways using the same data). 
This observation generalizes to any number of e-values, and holds for convex combinations or \emph{mixtures} that are not equally weighted. Of course, the e-values to mix and the weights for the mixture must be chosen without looking at the data; otherwise we are martingaling. Recently,~\cite{wasserman2020universal} used such averaging techniques to derandomize universal inference (discussed in Section~\ref{sec:ripr}). \cite{Vovk/Wang:2021} have shown that averaging is an admissible way of combining e-values (for a particular definition of admissibility) without further information about the e-values or their dependence structure. (And it is the only admissible symmetric method if we ignore the possibility of further mixing with the constant e-value 1.)

Test (super)martingales and e-processes can also be mixed, yielding mixture (super)martingales or e-processes. This \emph{method of mixtures} 
goes back to~\cite{Ville:1939}, \cite{Wald45}, and Robbins~\citep{darling1968some,robbins1970statistical,robbins1974expected}; see~\cite{howard2021time} for recent advances.

\subsection{Multiplying e-values}\label{subsec:mult-e}
Independent e-values can be combined by multiplication. If $B_1,\dots,B_n$ are independent e-values for $\cP$, then the product $B_1\cdots B_n$ is also an e-value for $\cP$.    
As we saw in Section \ref{subsec:test}, a product of dependent e-values can also be an e-value.  
If, for all $k$,  $B_k$ is an e-value for $\cP$ conditional on the values of $B_1, \ldots, B_{k-1}$, i.e.\ if
\begin{equation*}\label{eq:seqepre}
\E[B_k \mid B_1,\dots,B_{k-1}] \leq 1
\end{equation*}
for $k\ge1$, then $M_n = \prod_{k=1}^n B_k$ is an e-value for $\cP$. The sequence $(M_n)_{n\ge 0}$ is a supermartingale with respect to the filtration generated by the $B_k$.
In fact, in Section~\ref{subsec:test} we encountered, and in  Sections~\ref{sec:parametric} and~\ref{sec:nonparametric}
we again encounter, at each time $t$ a random variable $S_t$ which is a single-round e-variable conditional on the past, 
\begin{equation}\label{eq:seqe}
\E_P[S_t \mid \cF_{t-1}] \leq 1 \text{ for all $P \in \cP$}.
\end{equation}
If (\ref{eq:seqe}) holds for all $t$, then
$M_n = \prod_{t=1}^n S_t$ is an e-value for $\cP$. The sequence $(M_n)_{n\ge 0}$ is a supermartingale with respect to the (often richer) filtration $\cF$. 

\section{General principles and methodology}
\label{sec:principle}
As mathematical statisticians learned nearly a century ago from Jerzy Neyman and E.\ S.\ Pearson, the choice of a test of a null hypothesis should be guided by the alternative hypotheses that are considered plausible.  How should this work when we are using a test martingale, or more generally a test supermartingale or an e-process?

Intuitively, a good supermartingale should grow (get large) fast under the alternative so that we quickly build up evidence against the null as the sample size increases. So we want a test martingale or e-process with maximal expected rate of growth under the alternative.    

In this section, we first focus on testing a simple null $\cP=\{P\}$ against a simple alternative $\cQ=\{Q\}$ and use this case to develop our understanding of expected rate of growth (Section~\ref{sec:simple}). We then move to testing a simple null against a composite alternative (Section~\ref{sec:compositeH1}), and to the most difficult case, where even the null is composite (Section~\ref{sec:ripr}), introducing general methods for constructing e-processes---some based directly on growth rate optimality, some more indirectly. 

One danger we want to avoid throughout is an e-process becoming zero. Once this happens, the e-process can never become positive again, and thus it can never recognize later evidence against the null, no matter how strong. This can happen with positive probability under a particular alternative $Q$ only if the e-process's strategy for betting for $Q$ (i.e, the test martingale for $P$ designed to become large if $Q$ is correct; remember that the e-process is an infimum for such test martingales for the different alternatives) is sometimes allowed to bet all its money, thus risking bankruptcy.  We call this \emph{betting the farm}, and we insist on choosing e-processes that avoid it.

\subsection{Simple Null and Simple Alternative}
\label{sec:simple}
This is the case where we are testing a probability distribution $P$ against an alternative probability distribution $Q$.  As we saw in Section \ref{sec:evalues}, the likelihood ratio $dQ/dP$ is the natural test martingale in this case.  (This assumes that $Q$ is absolutely continuous with respect to $P$.)

What are the advantages of using this natural test martingale?  The most important advantage, perhaps, is that it has the greatest expected growth as measured using the \emph{expected logarithmic return}, a measure popularized by~\cite{Kelly56}. The name \emph{logarithmic return} is standard in finance and hence appropriate when we consider wealth processes.  When $E$ is a unit bet against $P$, $E$'s logarithmic return is simply $\log E$.  The fact that the expected logarithmic return $\E_Q(\log E)$ is maximized by $E:=dQ/dP$ can be obtained directly from Gibbs' inequality \cite[p.~413]{Shafer:2021}.
Because $M_t$ is a unit bet against $P$ whenever $M$ is a test martingale for $P$, it follows that the cumulative likelihood ratio
$E_t:=(dP/dQ)(X^t)$ maximizes 
\begin{align}\label{eq:gro}
    \E_Q( \log E_t)
\end{align}
and hence the growth rate $\E_Q(\log E_t)/t$
for each $t$.  By the same argument, $(dP/dQ)(X^\tau)$ maximizes $\E_Q( \log E_\tau)$ for every stopping time $\tau$.  \cite{GrunwaldHK19} called the requirement that $\E_Q( \log E_\tau)$ be maximized the \emph{GRO} criterion, for ``growth-rate optimality'' relative to $\tau$.  So we may summarize by saying that {\em in a simple vs.\ simple test, the likelihood ratio is GRO}. 

Why use the logarithm of $E$ rather than some other increasing function of $E$?  In finance, we average the logarithmic returns for successive time periods rather than the simple percentage returns in order to account for compounding.   \cite{Kelly56} pointed out that this compounding means that logarithmic returns add, and hence the law of large numbers applies, allowing us to gain some foresight about the medium to long run.  \cite{breiman1961optimal} showed that the logarithm has a number of other strong optimality properties, especially in iid\ settings where the wealth can be made to grow exponentially under the alternative, and this criterion maximizes the exponent. Using Wald's identity as Breiman used it, one can show that, in iid settings, the logarithm asymptotically (as $\alpha\to0$) minimizes expected time before $E$ reaches a desired threshold such as $1/\alpha$, independently of $\alpha$. It is also true that we will not ``bet the farm''\footnote{Advocates of Kelly betting in the stock market also use half Kelly and fractional Kelly strategies, which make sense when you are not confident about the alternative $Q$~\citep{Maclean2010Kelly}.} when we choose $E$ to maximize the expected logarithm, whereas this can happen if we maximize the expectation of $E$ itself or some polynomial function of $E$.  
See also \cite{Shafer:2021}, who compares expected logarithmic return to power in the Neyman-Pearson theory: both can be used to ask whether the alternatives for which a test is effective are plausible.

\subsection{Simple Null and Composite Alternative}\label{sec:compositeH1}

How do we find a good test martingale for $P$ when the alternative $\cQ$ is composite?  In general, we cannot maximize the expected growth rate under all the distributions in $\cQ$.  But we can look for an alternative $Q$ such that the test martingale defined by \eqref{eq:tmart} has a reasonably high expected growth rate under any distribution in $\cQ$ that fits the data $X_1,X_2,\dots$ reasonably well.  Because $X_1,X_2,\dots$ are revealed to us progressively, the natural procedure is to construct this $Q$ progressively.  On betting round $t$, we use the data so far, $x^{t-1}$, to choose the numerator $q(X_t\mid X^{t-1})$ in~\eqref{eq:unit}.  
Another way to view this
is to imagine the data being drawn from (or best explained by) some unknown $Q^* \in \cQ$ and --- since we do not know $Q^*$ --- to attempt to learn it from the data, at each round $t$ plugging in a $q(X_t\mid X^{t-1})$ that is an estimate of $q^*$ based on data $X^{t-1}$.

\subsubsection{The Plug-in Method.}

This is natural when $\cQ$ is a parametric model.  We use $x^{t-1}$ to estimate the parameters, and this gives us an estimate $\hat{Q}_t$ of the best fitting (or the `true') $Q^*$.  So our choice for $q(X_t\mid X^{t-1})$ is $\hat{q}_t(X_t\mid X^{t-1})$, where $\hat{q}_t$ is $\hat{Q}_t$'s density. Wald proposed, in passing and without any further analysis, this plug-in method \cite[Eq. 10.10]{Wald47}; it was subsequently analyzed by~\cite{robbins1974expected}, who connect it directly to the mixture method (introduced in the next subsection). 
Similar ideas were proposed independently by \cite{Dawid84} for prequential model validation and by \cite{Rissanen84} as a predictive version of Minimum Description Length  learning.
Recently, the plug-in method has been employed by~\cite{wasserman2020universal} in parametric models and~\cite{waudby2020confidence,waudby2020estimating} in nonparametric models.

We obtain a test martingale $M$ no matter how we estimate the $Q_t$. But we should not use maximum likelihood, at least when data are discrete, lest we end up betting the farm (the maximum likelihood estimator may assign probability $0$ to an outcome that may very well occur in the next round. Most of the authors just cited have found, however,  that it often suffices to slightly smooth the maximum likelihood estimator (often using a ``prior'') to avoid this problem, even in nonparametric settings.

\subsubsection{The Mixture Method.}
\label{subsubsubsec:mixmethod}

Another way to choose the $q(X_t\mid X^{t-1})$ in \eqref{eq:unit} is to average over the corresponding conditional distributions for the distributions in $\cQ$.  We can vary the weights with $t$ and with $X^{t-1}$.  This is the mixture method.
The mixture method is not a special case of the plug-in method, because the mixed probability distribution $Q$ we obtain may not be in $\cQ$ (this may happen if $\cQ$ is not ``fork-convex'', to use a concept introduced in \citealt{ramdas2022testing}). Since most models used in mathematical statistics are not fork-convex, $Q$ is rarely in $\cQ$.

\subsubsection{Bayes Factors.}\label{subsubsec:bfactor}

We can use a probability distribution on $\cQ$, say $R$, to define weights for a mixture martingale. We update $R$ on each round in the usual Bayesian way.  On round $t$, we use the update $R(\cdot\mid x^{t-1})$ in the averaging that produces $q(X_t\mid X^{t-1})$.

Not surprisingly, a simple calculation shows that the resulting unit bet $M_t=\prod_{i=1}^t B_i$ is equal to the Bayes factor defined by the distribution $R$ on $\cQ$ \citep{GrunwaldHK19}.  This Bayes factor is the ratio $q(X^t)/p(X^t)$, where $q$ is the density from mixing the distributions in $\cQ$ with $R$.  But for each $t$, the conditional probabilities given $X^{t-1}$ obtained by mixing with $R$ are the same as the conditional probabilities given $X^{t-1}$ obtained with $R(\cdot\mid X^{t-1})$.  Conditioning on the data so far commutes with averaging the distributions in the model.  

Bayes factors have been advocated by many statisticians as measures of evidence against a null when the alternative is composite \citep{berger1998bayes,Jeffreys61}.  E-values  measure evidence against the null in a different way.  Whereas a Bayes factor is used to multiply prior odds, an e-value is intuitively the outcome of a bet. Not surprisingly then, the correspondence between mixture test martingales (or e-values) and Bayes factors does not extend to composite nulls.

\subsubsection{Minimizing the Worst.}

If we do not  have a priori knowledge to guide us when determining   the `prior' distribution $R$ in the method of mixtures,
we may look for the distribution $R$  that minimizes the worst possible shortfall from this best growth rate.   This means that we measure the quality of test martingale $M$ stopped at time $\tau$ by
\begin{align}
    \label{eq:regrow}
\inf_{Q \in \cQ} \E_Q\left(\log M_\tau - \gamma \cdot \log M_\tau^Q\right).
\end{align}
where $\gamma =1$ and $M^Q$ is the GRO e-process so that $M^Q_{\tau} = (dP/dQ)(X^\tau)$, maximizing  (\ref{eq:gro}) with $t$ replaced by $\tau$.
We want this nonpositive quantity to be as large (as close to zero) as possible.  \cite{GrunwaldHK19} introduced this criterion and called it REGROW ({\em RElative GRowth Optimality in Worst case\/}).
A variation of this criterion, GROW, is obtained if we set $\gamma=0$. We then search for the $R$ that maximizes worst-case expected logarithmic return $\E_Q(\log M_{\tau})$ in an absolute rather than relative sense. In most applications, we are interested in REGROW rather than GROW optimality so we will focus on it below. Still, as will become clear in Section~\ref{sec:parametric}, in some applications GROW is more appropriate and in some applications they in fact coincide; \citep{GrunwaldHK19} discuss the differences in detail. 

Under slight regularity conditions, the e-variable $M_{\tau}$  that maximizes \eqref{eq:regrow} can be written as a Bayes factor defined relative to a specific prior $R_{\tau}$, where $R_{\tau}$ varies with $\tau$ and, in many cases, with $\gamma$. 
Still, in the cases considered by \cite{GrunwaldHK19} one can find priors $R$ that get us close to the maximum for all $\tau$. In some settings we do find a unique test martingale that maximizes \eqref{eq:regrow} for all $\tau$.  We will generalize \eqref{eq:regrow} to the case of composite nulls below, and we will find a unique e-process that maximizes the criterion for all $\tau$ for the group invariance tests of Section~\ref{sec:ttest}.

\cite{GrunwaldHK19} show that when $\cQ$ is an i.i.d.\ exponential family and $\cP = \{P\}$ is simple, the test martingale $M$ obtained from Jeffreys' prior is asymptotically REGROW: for every $\tau$ set equal to a large $t$, $M_t$ maximizes \eqref{eq:regrow} up to an $o(1)$ term 
among all e-variables that can be defined on $X^t$, linking growth optimality to description length (Section~\ref{subsubsec:mdl}).

\subsubsection{When Precise Alternatives are Hard to Come By.}
\label{subsec:grapa}
The previous development is based on the premise that we take the alternative $\cQ$ very seriously, as containing distributions from which the data may actually be sampled. More generally, we may simply have a set of distributions available which we think may predict the data reasonably well in terms of the log score, but which we would never believe to be `true', so that taking an expectation over them as in the GRO definition is not too meaningful. There is nothing that stops us from using such $\cQ$ in combination with e.g.\ the plug-in method; we will still have an e-process, we merely cannot claim its optimality any more. We return to this more `Fisherian' perspective on testing in Section~\ref{subsubsec:mdl}.  

Alternatively, 
\cite{waudby2020estimating}'s  GRAPA  (Growth Rate Adaptive to the Particular Alternative) method uses, on round $t$, the empirical distribution of $X^{t-1}$, possibly smoothed, for the numerator in \eqref{eq:unit}. In practice, this yields  efficient tests and (by inversion) confidence sets. 
GRAPA tries to mimic the growth rate that would be achieved by using a smoothed empirical distribution as the alternative (or its projection onto $\cQ$ when possible), while REGROW tries to match the $Q^*$-expected growth rate by learning $Q^*$ by the method of mixtures.

\subsection{Composite Null and Alternative}
\label{sec:ripr}

When the null $\cP$ and alternative $\cQ$ are both composite, we can usually handle them in a modular fashion.  The composite alternative can be handled as in the previous section, using the plug-in method or the method of mixtures.
Here we describe two relatively general ways of handling the composite null: universal inference (UI)---which always yields an e-process---and reverse information projection (RIPr)---which yields a sequence of e-variables that is sometimes an e-process. As mentioned in Section~\ref{subsec:test}, test (super)martingales for composite $\cP$ may not exist.  So we use the general concept of an e-process. 

UI and RIPr are not the only ways of handling composite nulls.  In Section~\ref{sec:nonparametric}, we will see many other test (super)martingales and e-processes for composite nulls. Most of these involve the method of mixtures applied directly to a collection of e-variables rather than to distributions in $\cQ$, as briefly introduced in Section~\ref{sec:alternativemixture} below.

\subsubsection{Universal Inference (UI).}

This method, introduced by \cite{wasserman2020universal} uses e-processes of the form 
\newcommand{\ui}{\ensuremath{\textsc{ui}}}
\begin{align}\label{eq:ui}
 \prode{t}{\ui} := 
 \frac{
 %
\bar{q}(X^t)}{\sup_{p\in\cP} p(X^t)} =
 \frac{\bar{q}(X^t)}{\hat{p}_{X^t}(X^t)},
\end{align}
where
$\bar{q}(X^t) := \prod_{i=1}^t \hat q_{X^{i-1}}(X_i)$, 
$\hat q_{X^{i-1}}$ is any distribution learnt from $X^{i-1}$, and $\hat{p}_{X^t}$ is the maximum likelihood estimator (MLE) under $\cP$, 
the final equality holding whenever the MLE is well-defined. 
Alternatively, we can use the method of mixtures and set  $\bar{q}(X^t) := \int \prod_{i=1}^t q(X_i  \mid X^{i-1}) dR(q)$, where $R$ is a distribution over $\cQ$.  In either case, as in the preceding subsection, the numerator is equal to $\bar{q}(x^t)$ for some alternative $\bar{Q}$ (usually not in $\cQ$), so that $\prode{t}{\ui}$ is the infimum of the family of test martingales $(\bar{q}(X^t)/p(X^t))_{p\in\cP}$ and hence an e-process by~\eqref{def:eprocess1}. 
The method is \emph{universal} because it does not require regularity assumptions or asymptotics and, importantly, is applicable in both parametric and nonparametric settings; see Section~\ref{sec:exch+LC} for an example of each.


We can think of $\prode{t}{\ui}$ as a middle ground between the non-Bayesian generalized likelihood ratio (MLE in both numerator and denominator) and the  Bayes factor for a composite null (mixtures in both numerator and denominator), neither of which leads to an e-process in general. By taking a supremum in the numerator, the generalized likelihood ratio exaggerates evidence for the alternative, requiring that this exaggeration be taken into account using the ratio's sampling distribution.  By including poorly fitting distributions in its mixture in the denominator, the Bayes factor may downplay evidence for the null.

\subsubsection{Reverse Information Projection (RIPr).}
This me\-thod, pioneered by \cite{GrunwaldHK19} (original arXiv version 2019) finds, for each stopping time $\tau$, an e-variable $\prode{\tau}{\textsc{ripr}}$ for $\cP$. 
To define $\prode{\tau}{\textsc{ripr}}$, we first choose a $\bar{Q}$ via the plug-in or mixture method, exactly like we did above for UI. 
%
Then we consider the set $\cW$ of all probability distributions on $\cP$, and for each $W\in\cW$, we denote by $P_W$ the distribution obtained by mixing the distributions in $\cP$ with $W$.

Extending results of \cite{Li99,LiB00},  \citet[Theorem 1]{GrunwaldHK19} show that, provided the infimum below is finite, for every $\tau \in {\cal T}$, there exists a unique measure $P^\tau$ for $X^\tau$ satisfying
\begin{equation}\label{eq:ripr}
D(\bar Q^\tau \|  P^\tau) = \inf_{W \in \cW} D( \bar Q^\tau \| P_W^\tau),
\end{equation}
where $D(\cdot \| \cdot)$ is Kullback-Leibler divergence and $\bar Q^\tau$ (resp. $P_W^\tau$) is $\bar Q^\tau$'s (resp. $P_W^\tau$'s) marginal for $X^\tau$. Further, $P^\tau$ has the following nontrivial property: defining
\[
  \prode{\tau}{\textsc{ripr}} := \bar q(X^\tau)/p^\tau(X^\tau),
 \]
where $p^\tau$ is the density of $P^\tau$, and $\bar q$ is the density of $\bar Q$,
$\prode{\tau}{\textsc{ripr}}$ is an e-variable; it is even the GRO e-variable relative to $\tau$, maximizing (\ref{eq:gro}) over all e-variables that can be written as a measurable function of $X^{\tau}$. 
$P^\tau$ is called the {\em reverse information projection} of $\bar Q$ onto $\cP$. In some cases (e.g.\ Section~\ref{sec:independence}, \ref{sec:ttest}) it is easy to calculate, in others (e.g.\ Cox regression \citep{TerschurePLG21}) it is not. In general, it is a sub-probability measure, i.e.\ $p^\tau$ may integrate to less than one; but in all cases of practical interest we have encountered so far, $P^\tau$ is a probability distribution. In particular, because of the convexity of KL divergence and of the set of mixtures of $\cP$, the infimum is often achieved by some $W\in\cW$ and then $P^\tau = P_W^\tau$.
The sequence $(\prode{t}{{\textsc{ripr}}})_{t\ge 1}$ is adapted to $\cF$. 

Sometimes (i.e.\ for some problem setting $\Pi,\cP,\cQ$ and choice of mixture or plug-in $\bar{Q}$)  it is an e-process; sometimes not.
If it is (as happens e.g.\ in the tests described in  Section~\ref{sec:ttest}), then  $(\prode{t}{{\textsc{ripr}}})_{t\ge 1}$  is an e-process relative to an appropriately chosen $\bar{Q}$. It then dominates UI when using the same mixture over $\cQ$:  $\prode{t}{\ui} \leq \prode{t}{\textsc{ripr}}$, since they have the same numerator, but UI maximizes  denominator likelihood. If it is not an e-process, then  variations of RIPr can often still be used to obtain one. For example, in parametric $k$-sample tests  (Section~\ref{sec:independence}),
$(\prode{t}{{\textsc{ripr}}})_{t\ge 1}$  is an e-process when for $\bar{Q}$ we take any fixed $Q \in \cQ$, but not when we take a plug-in or mixture $\bar{Q}$ that `learns'. To handle composite $\cQ$ we must then combine RIPr with the method of mixtures in a different way, as we now describe.  
\subsubsection{Mixing E-Processes.}\label{sec:alternativemixture}
In all nonparametric, and some parametric cases, a natural way to proceed is to first construct a parameterized collection of e-processes $\{\prode{}{\lambda}: \lambda \in \Lambda \}$. We then need to come up with a final e-process to use in practice. For this, we can use the method of mixtures again, but now by putting a distribution $R$ on the space $\Lambda$ and creating the new e-process  $\prode{}{R}$ with, for each $\tau$, $\prode{\tau}{R} := \int \prode{\tau}{\lambda} d R(\lambda)$, thus applying the method of mixtures directly to e-processes rather than to the alternative hypothesis $\cQ$, as we did above when introducing UI and RIPr. For example, in the two examples referred to above we have $\cQ = \{Q_{\theta}: \theta \in \Theta_1 \}$ a parametric set of distributions, and for each $\theta \in \Theta_1$, we define 
$
M_t^{\theta}:= (\prode{t}{\textsc{ripr}(Q_{\theta})})_{t \geq 1}
$
to be the sequence of RIPr e-variables relative to point alternative $Q_{\theta}$, which, as stated above, is an e-process.  
This is different from what we get if we try to use the method of mixtures directly to mix over $\cQ$, defining some $\bar{Q}$ that `learns' $Q_{\theta}$: in that case $(\prode{t}{\textsc{ripr}(\bar{Q})})_{t \geq 1}$ does not provide an e-process in the $k$-sample test of  Section~\ref{sec:independence} below. Indeed,  Section~\ref{sec:ttest} successfully mixes over $\cQ$, Sections~\ref{sec:independence} and all of Section~\ref{sec:nonparametric} mix over a collection of e-processes.
\section{Parametric Examples}\label{sec:parametric}
Examples of test martingales and e-processes for simple nulls abound in the Bayesian literature, since (a) every Bayes factor for a simple null also defines a test martingale. Further, (b) as pointed out by~\cite{darling1968some}, if $X_i$ are iid from $P$, and $P$ has a finite MGF, meaning $\Phi(\lambda):=\E_P[\exp(\lambda X_i)] < \infty$ for a given $\lambda>0$, then $\exp(\lambda \sum_{i\leq t}X_i) \Phi(\lambda)^{-t}$ forms a test martingale for $P$. Construction (a) amounts to positing a specific alternative (a mixture of elements of $\cQ$) and thus fits within both a Neyman-Pearson and Bayesian view on testing. In (b) there is no explicit alternative, and therefore it may perhaps be more in line with Fisher's view on testing. 

Below we emphasize examples with composite nulls, for which by now a plethora of e-variables have been designed. Some of these examples (logrank test, t-test, contingency tables) are implemented in the R package {\tt safestats} on CRAN \citep{ly2022safestats}.

\subsection{Reduction to a Simple Null}
One of the simplest practical ways to create e-processes with composite parametric nulls is to reduce the testing problem to one with a simple null, for example by coarsening  or conditioning on  a sufficient statistic of the data, so that the resulting marginal or conditional likelihood ratio has the same distribution under all $P \in \cP$, making this likelihood ratio an e-process. This strategy has already been applied in the early days of sequential testing with  Wald's sequential probability ratio test, which, while different from ours (Section~\ref{sec:seq}), is nevertheless based on a likelihood ratio between a simple null and a simple alternative. This led to sequential versions of the t-test (by marginalization, in 1950) and other group invariant testing problems as well as to a sequential  $2 \times 2$ contingency table test (by conditioning, in 1945). As reviewed in Section~\ref{sec:ttest}, it was recently discovered that for group invariant problems, such a reduction to a simple null leads  to e-processes that have both GROW and REGROW status, and are thus highly suitable for our anytime-valid context. As reviewed in Section~\ref{sec:conditioning}, the conditioning approach for the $2 \times 2$ table is suboptimal from a GRO perspective though. We review in Section~\ref{sec:independence} a general approach to e-processes for independence testing that does lead to an optimal e-process for the $2 \times 2$ table. But first we show that if $H_0$ and $H_1$ are separated, then optimal e-variables may sometimes simply be obtained by taking a likelihood ratio $q/p_0$ where $p_0$ represents a special element inside $\cP$.
\subsubsection{Monotone LR Families; Logrank Test.}
\label{sec:monotone}
\label{sec:logrankCHECK}
Fix a 1-dimensional regular exponential family given in terms of its densities $\{p_{\theta}: \theta \in \Theta \}$ with $\Theta$ representing either its mean- or its natural parameterization. 
Fix $\theta_0 < \theta_1$ and suppose that under the null, the $Y_i$ are i.i.d.\ $\sim P_{\theta}$ for $\theta \leq \theta_0$, whereas under the alternative, $\theta = \theta_1> \theta_0$. An easy calculation \cite[Example 4]{GrunwaldHK19} shows that the GRO e-variable for this problem is given by $E^{\theta}_t := \prod_{i \leq t} p_{\theta}(Y_i)/p_{\theta_0}(Y_i)$, i.e.\ it coincides with the likelihood ratio for a simple-vs.-simple testing problem. If the alternative is composite as well, i.e.\ $\theta\geq \theta_1$ for some $\theta_1 > \theta_0$, we have a choice: the GROW e-variable ((\ref{eq:regrow}) with $\gamma=0$) is given by setting $\theta= \theta_1$; alternatively, for a REGROW approach one can try to learn $\theta$ using the plug-in or mixture method as in Section~\ref{subsubsubsec:mixmethod}. 
The crucial property needed for these results to hold is a {\em monotone likelihood ratio property}, which does not only hold for exponential families but also, for example, if $p_{\theta}$ represents a t-distribution with a fixed degree of freedom and noncentrality parameter $\theta$. 

It also extends beyond the i.i.d.\ case: 
\cite{TerschurePLG21} use this idea to provide efficiently computable GROW e-variables and test martingales for the logrank test, a work-horse of medical statistics. The test martingale they derive can be extended to the more general setting of the Cox regression model with covariates, for which computationally efficient implementation remains a work in progress. For the simple logrank test without covariates they provide extensive simulations showing that whether GRO is preferable over REGROW is a subtle matter. 
{\color{black} Here we focus on this simple logrank test. One starts with $m_0$ subjects, partitioned into a {\em treatment\/}  $a$ and a {\em control\/} group $b$; for example, one wants to test a COVID vaccine; at time $0$ all $m_{a,0}$ subjects in the treatment group get the vaccine and all $m_{b,0} = m-m_{a,0}$ subjects in the control get a placebo. At every time $t=1,2,\ldots$ an {\em event\/} (e.g.\ onset of covid) happens, either to somebody in the treatment $(Y_i=a)$ or in the control $(Y_i=b)$ group. Under the null, treating patients does not yield any benefit, so that $P(Y_i=a \mid Y^{i-1}) = \theta m_{a,i}/(\theta m_{a,i}+ m_{b,i})$
for some $\theta \geq 1$ (if $\theta > 1$ then the treatment harms).  
Under the alternative, one assumes  that this {\em hazard ratio\/} $\theta$ satisfies $\theta \leq \theta_{\max}$ for some $\theta_{\max} < 1$.  Thus, at every $i$, one tests two separated sets of Bernoulli distributions, which can be addressed using the reduction to a null indicated above. }

\subsubsection{t-test, Regression, General Group Invariance.}
\label{sec:ttest}
Consider the following version of the t-test: according to the null, the $X_t$ are iid $\sim N(\delta_0 \sigma,\sigma)$ for some given effect size $\delta_0$; according to the alternative, they are iid $N(\delta_1\sigma ,\sigma)$ for effect size $\delta_1$. Under both null and alternative, the nuisance parameter  $\sigma$ is unknown, making the hypotheses composite. 
We coarsen the original process to $V_1, V_2, \ldots$, where $V_i := X_i /\lvert X_1\rvert$; of course, $\lvert V_1\rvert=1$. Under the null, $(V_t)_t$, by virtue of it being {\em scale-free\/} (it does not change if all data points are divided by a fixed constant), has a distribution $P_{\delta_0}$ that does not depend on the variance; similarly, under each distribution in the alternative, $(V_t)_t$  has the same marginal distribution $P_{\delta_1}$. So by considering $(V_t)_t$ instead of $(X_t)_t$ we reduce the problem to a simple-vs-simple test as in Section~\ref{sec:simple}, and the likelihood ratio $\prode{t}{} := p_{\delta_1}(V^t)/p_{\delta_0}(V^t)$ is a test martingale for the null relative to a coarsened filtration. 
Essentially the same likelihood ratio was proposed by \cite{Rushton50} for classical sequential testing.  \cite{Cox52} noted (using different terminology) that it can be rewritten as a Bayes factor applied to the original data under the improper right-Haar prior, $w(\sigma) = 1/\sigma$, i.e.
\begin{equation}\label{eq:haar}
\prode{t}{} = \frac{\int_{\sigma > 0} 
p_{\delta_1 \sigma,\sigma} (X^t)
w(\sigma) d \sigma}{\int_{\sigma > 0}  p_{\delta_0 \sigma,\sigma}(X^t) w(\sigma) d\sigma},
\end{equation}
where $p_{\mu,\sigma}$ denotes the density of a $N(\mu,\sigma)$ distribution.
\cite{Lai:1976} also noted the equality (\ref{eq:haar}) and first proposed to use $\prode{t}{}$ in an anytime-valid context, and even inverted the test to yield a closed-form confidence sequence for a Gaussian mean with unknown variance.

More recently \cite{perez2022estatistics} showed that, for all $\tau$, $\prode{\tau}{}$ is both the GROW and the REGROW e-variable:  among all e-variables for the given null it  maximizes
\begin{equation}
    \inf_{\sigma >0} \E_{P_{\delta_1 \sigma,\sigma}} \left(\log \prode{\tau}{} -  \gamma \log \prode{\tau}{\textsc{ripr}(P_{\delta_1 \sigma,\sigma})}
    \right), \nonumber 
\end{equation}
for every $\gamma \in {\mathbb R}$, including  $\gamma = 0$ (GROW) and $\gamma=1$ (REGROW), which is possible since $\E_{P_{\delta_1 \sigma,\sigma}} (\log \prode{\tau}{\textsc{ripr}(P_{\delta_1 \sigma,\sigma})})$ turns out to be constant in $\sigma$.  Note that this is the appropriate generalization of the (RE)GROW criterion 
((\ref{eq:regrow}) to the present case, with $\cQ$ set to $\{P_{\delta_1,\sigma,\sigma}: \sigma > 0 \}$ and 
since the null is now composite, $M^Q_{\tau}$ now set to 
$\prode{\tau}{\textsc{ripr}(P_{\delta_1 \sigma,\sigma})}$, which is now the GRO e-variable relative to $Q = P_{\delta_1 \sigma,\sigma}$ (Section~\ref{sec:ripr}).
This implies that the proposed e-process is in fact optimal in a strong sense.
They also demonstrate (RE)GROW optimality for the case (as in Section~\ref{sec:monotone} above)
that the defining constraint in either $H_0$ or $H_1$ or both is replaced by an inequality, i.e.\ $H_0$ expresses $\delta \leq \delta_0$ and/or $H_1$ expresses $\delta \geq \delta_1$, and for the case that, under $H_0$, $\delta=0$ while under $H_1$, $\delta$ is equipped with a prior distribution. In the latter case, if a heavy-tailed prior is used, the marginal of (\ref{eq:haar}) with respect to this prior coincides with the {\em Bayesian t-test\/} \citep{Jeffreys61,rouder-2009-bayes} which is thereby seen to have an e-process interpretation and to provide Type-I error safety. Thus, in this composite setting (and in contrast to the composite settings we discuss in the next sub-section), test martingales overlap with a specific popular type of Bayes factors.  
\cite{perez2022estatistics} also extend these insights to the general setting of $H_0$ and $H_1$ that share nuisance parameters expressing a group invariance: under a mild regularity condition, the (RE)GROW e-variable turns out to be equivalent to a likelihood ratio for a coarsening of the data, under which the null becomes simple.  
This includes scale invariance such as in the t-test above (if one divides all data points by the same constant, the e-variable for the t-test above remains invariant), but also location, rotation, and other invariances. By considering the affine group, one can handle linear regression problems with Gaussian errors as well.
\subsubsection{Conditional e-values.}
\label{sec:conditioning}
Suppose that, under the null, $Y \sim P_{\theta}$, $\theta \in \Theta$, and $Z$ is a sufficient statistic for $Y$, so that the density of $Y$ conditional on $Z$ is given by the same  $p(Y \mid Z)$, for all $\theta \in \Theta$. Then, for any conditional probability density  $q(Y \mid Z)$, 
we have by the same reasoning as used already in Section~\ref{sec:evalues}  that $S:= q(Y \mid Z)/p(Y \mid Z)$ is an e-variable conditional on $Z$, and hence also unconditionally, i.e.\ for all $\theta \in \Theta$, 
\begin{equation}\label{eq:conditioning}
\E_{Y \sim P_{\theta}} \left[ S \mid Z \right] = 1 \text{\ and so\ }
\E_{(Z,Y) \sim P_{\theta}} \left[ S  \right] = 1.
\end{equation}
Such e-variables were already used implicitly by \cite{Wald45}
in his approach to sequential independence testing in the $2 \times 2$ contingency table, in which the null is the Bernoulli model. The idea can be extended to more general exponential family nulls, but as  \cite{HaoGLLA23} show these usually do not have (RE)GRO(W) status. 

\subsection{Parametric Alternative, Nonparametric Null} 
\label{sec:independence}
A second broad class of e-variables arises when under the null, outcome $Y$ is independent of some observed covariate $X$ (or more generally, the null says that some measure of dependence takes on a value at most $\delta$); under the alternative, there is dependence (or more generally, the dependence is stronger than $\delta'$ for some $\delta'> \rho$). 
For simplicity assume $Y, Y_1, Y_2,\ldots $ are i.i.d.\ under both null and alternative, and, for now, assume that $X, X_1, X_2,\ldots $ are also i.i.d.\ (``random design'') according to some known distribution $R$, both under null and alternative. Then, with $\cP$ and $\cQ$ distributions for $(X,Y)$,  for all $P \in \cP$ , we have $P(X) = R(X)$, $P(Y\mid X) = P(Y)$,  whereas under all $Q \in \cQ$, $Q(X)= R(X)$ but $Q(Y \mid X) \neq Q(Y)$ with nonzero probability. Assuming that all $P \in \cP$ have density $p$, then for any conditional density $q(y \mid x)$, we define
\begin{equation}\label{eq:henzi}
s_q(X,Y) := \frac{q(Y \mid X)}{\int q(Y \mid x) dR(x)}.
\end{equation}
The random variable $S_q:= s_q(X,Y)$ is
trivially seen to be an e-variable for arbitrary $q$. \cite{GrunwaldHL22} show that 
$ \int q(Y \mid x) dR(x)$ is in fact the RIPr for the point alternative on $(X,Y)$ under which $X \sim R$ and $Y \mid X \sim Q \mid R$, and as a consequence, 
$S_q$ is  the GRO e-variable for $\cP$ as above and $\cQ= \{Q\}$, the simple alternative corresponding to density $q$. 
For such simple $\cQ$, it can be extended to a test martingale/e-process by 
setting $M_t := \prod_{i \leq t} s_q(X_i,Y_i)$. For composite $\cQ$, we can `learn' the right $q$ using the plug-in or mixture methods (Section~\ref{subsubsubsec:mixmethod}). Related ideas can be found in \cite{duan2020interactive} and \cite{Shaer2022}. 

This simple idea is surprisingly powerful: \cite{GrunwaldHL22} extend it to design e-processes for {\em conditional\/} independence testing in which there is a third random variable $Z$ and we test whether $P(Y \mid X, Z) = P(Y\mid Z)$ under the {\em model-X assumption} \citep{Candes2018}, which holds, for example, in randomized clinical trials. 
\cite{TurnerLG21,TurnerG22} provide another, quite different extension, without $Z_i$, in which $X_i$ is not random. Rather, it may take on a finite number, say $k$,  of values, and we observe $k$ {\em data streams}. To illustrate, consider the case $k=2$ with  $X_i \in \{a, b\}$ binary, $X_i= a$ denoting that $Y_i$ is an outcome in the `treatment' group and $X_i=b$ that $Y_i$ is in the `control' group. We then observe streams 
$Y_{1,a}, Y_{2,a}, \ldots$ and $Y_{1,b}, Y_{2,b}, \ldots$. 
This is a sequential two-sample test (revisited in Section~\ref{sec:2st}): the null now expresses that both streams are i.i.d.\ according to the same $P$, under the alternative, the distributions differ. 
Although now the $X_i$ are not i.i.d.\ any more, we can `flatten' the e-variables (\ref{eq:henzi}) so that they can still be used in this context. \cite{TurnerLG21} show that variations of the resulting  e-variables enjoy (RE)GROW optimality properties. \cite{TurnerG22} further extend the setting to the case that the null does not express independence but rather that the strength of dependence is bounded by some user-supplied effect size $\delta$ (such as difference in mean or, if $Y_i \in \{0,1\}$, the odds-ratio), and use this to develop anytime-valid confidence sequences. 
Finally, \cite{TurnerG23} allow addition of {\em strata\/} $Z$ thus providing a safe Cochran-Mantel-Haenszel test. 
\commentout{
Similar sequential $k$-sample tests, but with the smaller null that $P(Y \mid X=a) = P(Y \mid X=b) = P$ where $P$ is restricted to be a member of a given exponential family,  can also be implemented using products of conditional e-variables as in (\ref{eq:conditioning}).
\cite{HaoGLLA23} extensively compares such conditional e-variables to the e-variables based on (\ref{eq:henzi}) and finds that the differences in growth-rate under a specific exponential family alternative is usually quite small. 
}

\section{Nonparametric Examples}\label{sec:nonparametric}

We present case studies that illustrate how one builds  
\begin{itemize}
    \item[(A)] test martingales for composite nonparametric $\cP$ despite there being no common reference measure, 
    \item[(B)] test \emph{super}martingales when no martingales exist, 
    \item[(C)] e-processes when no supermartingales exist, and 
    \item[(D)] confidence sequences (CSs) for functionals using \emph{reversed} submartingales instead of test martingales. 
\end{itemize}
They are presented in an order that aids readability due to complexity of concepts, rather than (A) to (D) above. Many of the CSs have been implemented in C++ in the package \texttt{confseq}\footnote{https://github.com/gostevehoward/confseq}, with Python and R interfaces.

\subsection{Estimating Sub-Gaussian Means (Case B)}
\label{sec:subgaussian}
A distribution $P$ for real-valued $X_1,X_2,\dots$ is \emph{sub-Gaussian} with parameter $\sigma>0$ if 
\[
  \forall\lambda,i:
  \E_P[\exp(\lambda(X_i - \mu_i)) \mid \cF_{i-1}]\leq \exp(\lambda^2\sigma^2/2),
\]
where $\mu_i := \E_P[X_i \mid \cF_{i-1}]$.
Fixing $\sigma$, let $\cG^\mu$ be the set of sub-Gaussian distributions with parameter $\sigma$ and $\mu_i=\mu$ for all $i$; 
this is a nonparametric generalization of Gaussianity.
\cite{darling1968some} constructed CSs for $\mu$. 
They observed that for any $\lambda\in\R$,
\begin{equation}\label{eq:sub-Gaussian-superMG}
M^\mu_t(\lambda) := \exp\left(\lambda \sum_{i \leq t} (X_i-\mu) - \tfrac{\lambda^2}{2} \sigma^2 t \right)
\end{equation}
is a test supermartingale for $\cG^\mu$. Setting $Y_t := \sum_{i \leq t}X_i$, and choosing a centered Gaussian with variance $\rho^2$ as a mixing distribution $F$, they define
\[
M^\mu_t := \int M^\mu_t(\lambda) dF(\lambda) = \frac{\exp(\tfrac{\rho^2(Y_t -t\mu)^2}{2(t\sigma^2\rho^2+1)})}{\sqrt{t\rho^2\sigma^2+1}},
\]
which is also a test supermartingale for $\cG^\mu$.  It grows exponentially fast under any $P\in \cG^\theta$ for any $\theta\neq \mu$, and it grows faster for $\theta$ farther from $\mu$; thus, the evidence automatically adapts to the difficulty of the testing problem. This type of adaptivity is a commonly observed benefit of the method of mixtures, if appropriately employed. 

Using the inversion~\eqref{eq:general-CS-by-inversion-2}, we find that 
\begin{equation}\label{eq:dr}
\frac{Y_t}{t} \pm \sigma \sqrt{\frac{(t\rho^2+1)}{t^2\rho^2}\log((t\rho^2+1)/\alpha^2)}
\end{equation}
is a CS for $\mu$. If $\mu_i$ differs for each $i$, an identical argument shows that~\eqref{eq:dr} is a CS for the running mean $\sum_{i=1}^t \mu_i/t$. The general scaling of $\sigma t^{-1/2} \sqrt{\log t + \log\alpha^{-1}}$ is expected when using mixture distributions that are continuous around the origin \citep[Prop. 2]{howard2021time}. The $\sqrt{\log t}$ can be changed to $\sqrt{\log\log t}$ at the expense of other constants using mixture distributions that are unbounded at the origin; see \citet[Eq.~(1)]{howard2021time}. 

\cite{waudby2020estimating} observed that 
\[
  \exp\left(\sum_{i \leq t}
  \left(\lambda_i (X_i-\mu) - \tfrac{\lambda_i^2}{2} \sigma^2 \right) \right)
\]
is also a test supermartingale for $\cG^\mu$, whenever $\lambda_i$ is predictable. They invert  test supermartingales of this form to yield statistically efficient CSs.

\subsection{Heavy-Tailed, Robust Mean Estimation (Case~B)}
\label{sec:heavy}
For a fixed $\sigma>0$, let $\cV^\mu \supset \cG^\mu$ be the set of distributions on $\R^\infty$ that yield observations with conditional mean $\mu$ and conditional variance bounded above by $\sigma^2$.  
Inspired by~\cite{catoni2012challenging}, \cite{wang2022catoni} prove that
\[
L^\mu_t := \exp\left( \sum_{i \leq t}\varphi( \lambda (X_i-\mu)) -  \tfrac{\lambda^2}{2} \sigma^2 t \right)
\]
is a test supermartingale for $\cV^\mu$, where $\varphi(x)$ equals $\log(1+x+x^2/2)$ if $x\geq0$ and $-\log(1-x+x^2/2)$ if $x<0$. They then use the plug-in and inversion techniques discussed in the previous cases to derive a CS for $\mu$.
Somewhat surprisingly, these CSs for $\mu$ (assuming $P \in \cV^\mu$) appear, visually, almost identical to the sub-Gaussian CSs one gets when assuming $P \in \cG^\mu$. In other words, for mean estimation, the sub-Gaussian assumption can be relaxed with almost no practical consequence. 
The authors also derive extensions for the case when the $p$-th moment is finite ($p>1$).

There are other known test supermartingales for this setting, for example by~\cite{dubins1965tchebycheff} and~\citet[Proposition 12]{delyon2009exponential} but experiments show that these to be significantly less powerful. 

\cite{wang2023huber} later extended these ideas to derive ``Huber-robust'' test supermartingales and CSs that can handle adversarial corruptions to  heavy-tailed data.

\subsection{Variance-Adaptive Estimation of Bounded Means (Case A)}
\label{sec:bounded}

In the previous two examples, the sub-Gaussian parameter or the variance bound $\sigma$ must be provided (or an upper bound must be guessed) in advance by the statistician, since it is provably impossible to learn $\sigma$ from the data itself. 
For the subclass of bounded random variables, however, variance-adaptive mean estimation is feasible.

Let $\cB^\mu$ denote the set of distributions $P$ on $[0,1]^\infty$ such that $\E_P[X_i\mid\cF_{i-1}]=\mu$.
\cite{howard2021time} prove that for any $\lambda \in [-1,1]$, and any predictable $\hat \mu_i  \in \cF_{i-1}$,
\begin{equation}
\label{eq:subexponential-superMG}
N^\mu_t(\lambda) := 
\exp\left(\lambda \sum_{i \leq t} (X_i-\mu) - \psi(\lambda) \sum_{i\leq t} (X_i - \hat \mu_i)^2 \right),
\end{equation}
where $\psi(\lambda):= -\log(1-\lambda)-\lambda$, is a test supermartingale for $\cB^\mu$. Because $\psi$ is the logarithm of the moment generating function (MGF) of a centered unit-rate exponential distribution, we call $N^\mu$ a subexponential supermartingale. As $\lambda\to0$, $\psi(\lambda)$ behaves like the Gaussian log-MGF $\lambda^2/2$, but unlike the sub-Gaussian supermartingale in~\eqref{eq:sub-Gaussian-superMG}, we can  employ a fully empirical variance term in~\eqref{eq:subexponential-superMG}. This generalizes a result of~\cite{fan2015exponential}, who effectively proved the same claim with $\hat \mu_i:=0$ for all $i$. The extension to predictable $\hat \mu_i$, obtained using some tricky algebra, is useful in lowering the empirical variance. Mixing over $\lambda$ with a (conjugate) gamma distribution leads to a closed form mixture supermartingale; see~\cite{howard2021time} for the resulting CS.


\cite{waudby2020estimating} note that for predictable $\lambda_i \in \cF_{i-1}$, 
\[
N^\mu_t := 
\exp\left( \sum_{i \leq t} \lambda_i (X_i-\mu) -  \sum_{i\leq t} \psi(\lambda_i)(X_i - \hat \mu_i)^2 \right)
\]
is also a subexponential ``plug-in'' test supermartingale that can be tuned to closely mimic the earlier mixture supermartingale. This means that
\[
\frac{\sum_{i \leq t} \lambda_i X_i}{\sum_{i \leq t} \lambda_i} \pm \frac{\log(2/\alpha)+\sum_{i\leq t} \psi(\lambda_i) (X_i - \hat \mu_i)^2}{\sum_{i \leq t} \lambda_i}
\] 
is a $(1-\alpha)$-CS for $\mu$.  

The preceding techniques are interesting because they can be used even when the observations are not bounded. But for the bounded model $\cup_{\mu\in\R}\cB^\mu$, the most statistically powerful way to derive a CS for $\mu$ is to use plug-in test martingales for $\cB^\mu$ of the form 
\begin{equation}\label{eq:product-capital-bounded}
K^\mu_t := \prod_{i=1}^t (1+\lambda^\mu_i(X_i-\mu)),
\end{equation}
where $\lambda^\mu_i$ is a predictable process indexed by $\mu$. As before, $C_t:=\{\mu: K^\mu_t < 1/\alpha\}$ is a $(1-\alpha)$-CS for $\mu$. The $\lambda^\mu_i$  are naturally interpreted as bets on the $X_i$; they must be predictable because a bet on $X_i$ must be made before seeing $X_i$. 
This idea was suggested  by~\cite{hendriks2018test}, and was independently proposed and studied in more depth by \cite{waudby2020estimating}, who derive betting strategies that are adaptive to the underlying distribution $P$, in particular to its mean and variance, establishing connections to the Chernoff method, empirical and dual likelihood, and other parts of the literature.  
Followup work by~\cite{orabona2021tight} derives other betting strategies via connections to Thomas Cover's universal portfolios. 

\subsection{Testing Symmetry (Case A)} 
\label{sec:symm}

Let $\cP$ be the set of distributions on $\R^\infty$ such that $X_t$ and $-X_t$ have the same distribution given $\cF_{t-1}$, for every $t\geq 1$.
Extending an older result by \cite{efron1969student}, 
\citet[Lemma 6.1]{delapena:1999} establishes that for any $\lambda \in \R$,  
\begin{equation}
\label{eq:symm-supermg}
R_t(\lambda) := \exp\left(\lambda \sum_{i \leq t}X_i - \tfrac{\lambda^2}{2}\sum_{i \leq t} X_i^2\right)
\end{equation}
is a test supermartingale for $\cP$. Notice again the fully empirical variance term, as in~\eqref{eq:subexponential-superMG}; these are also called self-normalized processes. As before, one can mix over $\lambda$ or use the plug-in technique described earlier.

Recently, \cite{ramdas2020admissible}
proved that $R_t$ is inadmissible for testing symmetry by constructing a test \emph{martingale} $R^{o}_t$ for $\cP$ that is always at least as large as $R_t$, and typically larger. 
In fact,  $M$ is a test martingale for $\cP$ (and is admissible) if and only if the unit bets $B_t$ at time $t$ in~\eqref{eq:unit} are nonnegative and predictable, and $B_t-1$ is an odd function of $X_t$. The unit bet underlying~\eqref{eq:symm-supermg} takes the form~$g(x):=\exp(\lambda x - \lambda^2 x^2/2)$. Since $g(x)-1$ is not an odd function, mirroring it around one can improve it: using unit bets $\widetilde{g}(x) := g(x)1_{x \geq 0} + (2-g(-x))1_{x\leq 0}$ yields a strictly better (and admissible) test martingale $R^o_t$.

\subsection{Testing Exchangeability and Log-Concavity (Case C)}
\label{sec:exch+LC}

In the previous example, $\cP$ is a  rich, nonparametric class of distributions (discrete and continuous, light and heavy tailed, etc.) with no common dominating measure. Being able to find a single (nonconstant) process that is simultaneously a test martingale for every $P$ in $\cP$ is quite atypical. (The same atypical situation also occurred with $\cB^\mu$ in the bounded case.) For example, there is no nontrivial test martingale for $\cG^\mu$, the sub-Gaussian class discussed earlier; nevertheless, we did exhibit a test supermartingale. It turns out that even this is atypical: a rather special structure is required for even a (nontrivial) test supermartingale to exist. 

\cite{ramdas2022testing} study the seemingly simple problem of testing if a binary sequence is exchangeable, and find that \emph{no nontrivial test supermartingale exists} (in the original filtration), but they exhibit a nontrivial and powerful e-process based on universal inference.  

Remarkably, \cite{vovk2021testing} demonstrates that by \emph{shrinking the filtration} to include only conformal p-values, it is once again possible to design nontrivial test martingales, even though none exist in the richer data filtration. Vovk's method works for general observation spaces, but in the binary case, experiments by~\cite{vovk2021conformal} demonstrate that it is not as powerful as the aforementioned e-process.

Another relevant example is that of testing log-concavity. Let $\cL_d$ denote the set of distributions $P$ on $\R^d$ with Lebesgue densities $p$ such that $\log p(x)$ is concave in $x$. $\cL_d$  is a nonparametric class that contains all Gaussian, logistic, exponential and Laplace distributions, as well as uniform distributions on any convex set. \cite{gangrade2023sequential} prove that there is no test supermartingale for $\cL_d$, but the universal inference approach yields a powerful e-process for $\cL_d$.

Of course, there are problems for which even no nontrivial e-process exists and testing those nulls is futile; see~\cite{ruf2022composite} for examples.

\subsection{Estimating Convex Functionals and Divergences by Reversing Time (Case D)} 

Consider  a set of probability distributions $\Pi$ that is closed under convex combinations. 
A functional $\phi: \Pi \mapsto \R_{\geq 0}$ is called convex if $\phi(a P + (1-a)Q) \leq a\phi(P) + (1-a)\phi(Q)$ for any $P,Q\in\Pi$ and $a\in[0,1]$. Classic examples are the entropy and the mean. Similarly, a divergence $D: \Pi \times \Pi \mapsto \R_{\geq 0}$ is called convex if $D(a P + (1-a)P', aQ + (1-a)Q') \leq aD(P,Q) + (1-a)D(P',Q')$. Examples include the total variation distance, Kullback-Leibler divergence, kernel maximum mean discrepancy, Kolmogorov-Smirnov distance, Wasserstein distance or any integral probability metric or f-divergence. 

Suppose $X_1,X_2,\dots,X_t,\dots \sim P$ and let $P_t$ denote the empirical distribution of $X^t$. The exchangeable filtration $\mathcal E_t$ is the \emph{decreasing} filtration given by $\mathcal E_t=\sigma(P_t,X_{t+1},X_{t+2}\dots)$; in words: $X_{t+1},X_{t+2}\dots$ are known perfectly, but the order of $X_1,X_2,\dots,X_t$ is forgotten.
\cite{manole2021sequential} derive a curious property: for any convex functional $\phi$, the process $(\phi(P_t))_{t \geq 1}$ is a \emph{reverse} submartingale with respect to the exchangeable filtration. (An analogous statement also applies to divergences.)

Recall that a reverse submartingale is a submartingale when time is reversed and the process is viewed from time $\infty$ to zero.  Reverse submartingales behave somewhat like forward supermartingales: their expectations are decreasing as time increases. Nonnegative reverse submartingales behave like test supermartingales in that there exists a \emph{reverse} Ville's inequality, with an identical statement to the forward Ville's inequality.  \cite{manole2021sequential} use this to derive confidence sequences for (say) the entropy of a distribution, as well as for divergences between pairs of distributions in quite some generality. The same technique also allows the authors to derive the first tight CSs (in their dependence on sample size, dimension, etc.) for suprema of Gaussian processes, Rademacher complexities, U-statistics, quantile functions, and several other interesting objects.

A game-theoretic interpretation of nonnegative reverse submartingales remains unknown, as does a game-theoretic derivation of the above CSs.

 \subsection{Sequential Change Detection} \label{sec:change}
 
 On observing a stream of data, the problem of sequential change detection can be seen as an extension of sequential testing: either all the data is from some $P\in\cP$, or at some time $\nu$, it switches from $P$ to some $Q \in \cQ$. If there is indeed a change, we would like to stop as quickly as possible and proclaim a change; if this happens, call the time at which a change is proclaimed as $\tau^*$. Measures of the performance of a change detection procedure include average run length (ARL, also called frequency of false alarms) and average detection delay. These are respectively defined as $\inf_{P \in \cP}\E_P[\tau^*]$, and $\sup_{P \in \cP,\nu > 0,Q \in \cQ}\E_{P,\nu,Q}[\tau^* - \nu \mid \tau^* > \nu]$, where the subscript $P,\nu,Q$ means that the data come from $P$ up through time $\nu$ and from $Q$ after $\nu$. We would like the former to be as large as possible with the latter being as small as possible. 

Extending \cite{volkhonskiy2017inductive}, who use conformal test martingales to detect deviations from exchangeability of the $X_i$'s, 
\cite{shin2022detectors} describe a general nonparametric game-theoretic framework for change detection. They define an \emph{e-detector} for $\cP$ to be a nonnegative process $M$ such that 
\[
\E_P[M_\tau] \leq \E_P[\tau]
\text{ for every $\tau \in \Tau$ and $P \in \cP$}.
\]  
(Like an e-process, the definition only depends on $\cP$ but we measure its quality relative to a post-change class $\cQ$). The authors show that if one can construct an e-process for $\cP$, then one can define an e-detector for $\cP$ by summing e-processes started at consecutive times. Formally, $M_t := \sum_{i \leq t} A_i$ is an e-detector for $\cP$, where $A_i$ is an e-process for $\cP$ that depends only on $X_i,X_{i+1},\dots$. Game-theoretically, $M_t$ is the wealth of a gambler who injects an extra dollar into the game at each time, and uses it to bet against $\cP$.

The above definition and construction may appear mysterious, but yields methods with nontrivial properties. First, defining $\tau^* := \inf\{t \geq 1: M_\tau \geq 1/\alpha\}$, one can prove that the ARL is at least $1/\alpha$. Second, in certain parametric settings, if there is indeed a changepoint, then one can design e-detectors such that the detection delay is near-optimal in a particular sense: even if $P,Q$ were known in advance, the best possible detection delay for any method with ARL at least $1/\alpha$ scales like $\log(1/\alpha)/D(P||Q)$, where (as before) $D$ is the Kullback-Leibler divergence. An e-detector based on likelihood ratios recovers the famous Shiryaev-Roberts statistic and can adaptively achieve this optimal scaling (up to lower order terms) without knowing $P,Q$, by employing new mixture and plug-in approaches. Third, e-detectors can be built for many nonparametric $\cP$ using (for example) the e-processes constructed earlier in this section. For many such nonparametric problems, e-detectors provide, as far as we know, the first change detection procedures with provable ARL control.

E-detectors only work when $\cP$ and $\cQ$ are prespecified and non-intersecting (eg: we wish to detect mean changes from $\leq 0$ to $>0$).
 \cite{shekhar2023sequential} develop a complementary framework for change detection when only such information is not available (eg: we wish to detect mean changes from anything to anything else) using a new notion of a ``backward confidence sequence'' (BCS). Here, one constructs two $(1-\alpha)$-confidence sequences --- one forwards in time, and one backwards in time (i.e.\ forwards in time after reversing time) --- and declares a changepoint the first time that they do not intersect. One can show that the ARL is controlled nonasymptotically, and in parametric settings one achieves an optimal detection delay up to constants. Due the plethora of CSs developed for nonparametric settings, plugging them into the BCS framework yields the first sequential change detection method for many nonparametric problems.

\subsection{Time-Uniform Central Limit Theory and \emph{Asymptotic} Confidence Sequences}

The average treatment effect (ATE) is arguably the most popular estimand in causal inference, and one may ask if it is possible to estimate it sequentially. For brevity, we focus here on the observational setting, where finite-sample inference is not possible due to unknown biases caused by confounding, but under suitable assumptions it is possible to design ``doubly-robust'' estimators for the ATE that have (nonsequential) asymptotic coverage guarantees. For a sequential analog, a suitable generalization of the concept of confidence sequences is required, because (by definition) CSs have finite-sample validity guarantees that are already impossible in nonsequential settings.

With such goals in mind, \cite{waudby2021time} define ``asymptotic confidence sequences'', which may sound paradoxical at first. They mirror an analogous definition of asymptotic CIs. Informally, a sequence of (measurable) sets $(C_t)_{t \geq 1}$ is called an asymptotic CS if there exists some unknown nonasymptotic CS $(D_t)_{t \geq 1}$, such that the measure of the symmetric difference between $C_t$ and $D_t$ almost surely vanishes faster than a $\sqrt{\log\log t/t}$ rate. 
\cite{waudby2021time} then derive a universality result: informally, as long as the data have more than two moments (an almost necessary condition for inference), a universal asymptotic CS is given by~\eqref{eq:dr} but with $\sigma$ replaced by an empirical variance $\hat \sigma_t$. This yields a time-uniform analog of the central limit theorem (CLT), and is established using certain strong approximation theorems for Brownian motion. The authors then construct doubly-robust asymptotic CSs for the ATE, yielding anytime versions of the corresponding CIs.

Asymptotic CSs can be used in a variety of other settings where CLT-based CIs are the norm in the offline setting. These include M-estimation and other semiparametric and nonparametric functional estimation problems; see also~\cite{pace2020likelihood} and~\cite{johari2022always}.


In complementary work,  \cite{duan2022interactive} develop test martingale versions of rank based tests like the Wilcoxon, Kruskal-Wallis, and Friedman tests. Batch versions of these tests are commonly used for testing the strong global null (of no treatment effect) in a randomized experiment with covariates.

\subsection{Other Nonparametric Problems}

\subsubsection{Compendium of Exponential Supermartingales.}
In previous sections, we have encountered several nonparametric test supermartingales of the form
\[
\exp(\lambda \cdot (\text{sum}_t) - \psi(\lambda)\cdot(\text{variance}_t) ),
\]
where $\text{sum}_t$ is the sum of the observed random variables, and $\text{variance}_t$ captures their cumulative variance; recall~\eqref{eq:sub-Gaussian-superMG},~\eqref{eq:subexponential-superMG},~\eqref{eq:symm-supermg} for example.
Likelihood ratios for in exponential families also take an identical form, where $\text{sum}_t$ adds the sufficient statistics. In a very concrete sense, we have been generalizing the likelihood ratio to composite and nonparametric problems, and even to cases where there is no dominating measure. Just as likelihood ratios are fundamental objects for parametric inference, test (super)martingales and e-processes are fundamental objects for nonparametric inference.
\cite{howard2020time} summarize a large literature on test supermartingales and e-processes of the above exponential form, in discrete and continuous time, for scalar-, vector- and matrix-valued observations, and under a variety of nonparametric conditions. We have mentioned only some examples.

\subsubsection{Quantiles.} 
\cite{howard2022sequential} derive confidence sequences based on iid\ data for any prespecified quantile of an unknown probability distribution, improving on those derived by \cite{Darling/Robbins:1967b}. They also derive CSs for the entire cumulative distribution function (or quantile function) of any arbitrary univariate random variable, proving a time-uniform extension of the famous Dvoretzky-Kiefer-Wolfowitz inequality.

\subsubsection{Two-Sample (and Independence) Testing.}\label{sec:2st}
Here, we observe two samples and want to know if they have the same distribution, making no other assumptions. This is one of the best studied problems in statistics.  Nonparametric
methods in the offline setting include the univariate Kolmogorov-Smirnov test and the multivariate kernel maximum mean discrepancy, amongst many others. 
\commentout{
However, the literature on sequential nonparametric two-sample testing appears sparse: \cite{balsubramani2015sequential} and~\cite{lheritier2018sequential} do use test (super)martingales, but the methods they proposed are not that powerful in practice. 
\cite{shekhar2021game} extend these ideas in a relatively general game-theoretic framework that provides the first sequential analog of large classes of offline nonparametric tests (even for non-i.i.d.\ data), and these perform excellently in practice. The evidence grows slowly for hard problems (when the two distributions are different but very similar) and quickly for easy ones (when the two distributions are very different), and it can be monitored and stopped adaptively. This is a major advantage over offline tests when the problem difficulty is not known in advance. Extending the above, the first sequential nonparametric independence testing framework was developed in~\cite{Podkopaev22skit}, which allows random variables to lie in general spaces and also handles non-i.i.d.\ settings.
}

As to sequential two-sample tests, we should distinguish between three approaches. 
We already discussed the two-sample  tests from  Turner et al.  (\citeyear{TurnerG22,TurnerG23}) in which the alternative is parametric in Section~\ref{sec:independence}. Although we called these tests `parametric', their null is better viewed as nonparametric. Second, in some works we can use arbitrary (e.g.\ deep learning based) sequential predictors that attempt to predict, given an outcome, from which of the two samples it was taken.  
These include \cite{lheritier2018sequential} (although the likelihood ratio process they use technically is not an e-process) 
and \cite{Pandeva22} (which may be re-interpreted as a  modification of the above process based on UI, so that it does become an e-process); 

Finally, one may directly attempt to obtain sequential analogs of large classes of existing offline nonparametric tests. 
\cite{shekhar2021game} successfully pursue this direction (even for non-i.i.d.\ data) in  a general game-theoretic framework, and their e-processes perform excellently in practice. The evidence grows slowly for hard problems (when the two distributions are different but very similar) and quickly for easy ones (when the two distributions are very different), and it can be monitored and stopped adaptively. This is a major advantage over offline tests when the problem difficulty is not known in advance. Extending the above, the first sequential nonparametric independence testing framework was developed in~\cite{Podkopaev22skit}, which allows random variables to lie in general spaces and also handles non-i.i.d.\ settings.

\subsubsection{Sampling Without Replacement (WoR).}
Another classical problem is that of estimating a mean when sampling WoR. Here we have a bag of $N$ numbers $\{x_1,\dots,x_N\}$, say all in the range $[0,1]$, and we wish to estimate their average $\mu:=\sum_{i \leq N} x_i / N$, or (say) to test if it is at most a half. The randomness arises from the WoR sampling process. \cite{waudby2020confidence} construct powerful plug-in test supermartingales for testing such hypotheses (of the empirical Bernstein flavor in~\eqref{eq:subexponential-superMG}), and invert them to construct CSs. \cite{waudby2020estimating} designed more powerful test martingales of the form~\eqref{eq:product-capital-bounded} and the resulting confidence sequences are the tightest known so far. These were then applied quite successfully towards election auditing by~\cite{waudby2021rilacs} and more recently by~\cite{spertus2022sweeter}.


\section{Multiple Hypothesis Testing}

\subsection{Global Null Testing and Meta-Analysis}\label{sec:meta}
Based on test martingales, \citet{SchureG22} propose {\em ALL-IN\/} ({\em A\/}ny time {\em L\/}ive and {\em L\/}eading {\em IN\/}terim) meta-analysis. This meta-analysis can be updated {\em any time}, even after each new observation,  while retaining type-I error guarantees.  It is {\em live\/}: no need to specify in advance the times when you will look and reanalyze.  And it can be the {\em leading\/} source of information for deciding whether individual studies should be initiated, stopped early, or expanded.  

These authors illustrate the method for clinical trials involving time-to-event data, using a Gaussian approximation to Section~\ref{sec:monotone}'s logrank test. 
Consider the case where each study tests the null hypothesis that some effect size $\delta$ (measuring, say, the efficacy of a medical treatment) is $0$; extensions to CIs are possible via inversion. In the simplest case, the evidence for the $i$-th study is measured by a unit bet $\seqe{(i)}{}$; the null is always the same (a ``global null''), but the alternative may change. For example, if the first study is based on the mixture method of Section~\ref{sec:compositeH1}, the mixing distribution for later studies might be updated using the outcomes of the studies so far or changed because the next study samples from a different population.
The unit bets $\seqe{(1)}{}, \seqe{(2)}{}, \ldots$ generated this way can be multiplied, so that the process $\prode{}{}$ with $\prode{(j)}{} := \prod_{i=1}^j \seqe{(i)}{}$ is a test martingale at the ``meta-level'', with individual outcomes replaced by entire studies. 
We can always keep initiating and adding new studies as we want at the time, deciding whether do so and choosing the unit bet for any new study in light of the outcomes of the previous studies.%

In the terminology of \cite{GrunwaldHK19}, the method is safe under {\em optional continuation}. This is true when  the study-level e-variables $\seqe{(j)}{}$
%
are produced by Section~\ref{sec:ripr}'s RIPr,  even when the RIPr does not give an e-process at the individual outcome level.  The method is  more flexible though when each study $j$ is associated with an e-{\em process}. As \cite{SchureG22} show, it is then possible to interleave the studies---one may first observe some outcomes from study 1, then some from study 4, then some from study 1 again, etc., tracking the cumulative product of the e-variables resulting from each batch. Again, one can decide at any time to stop an individual study, initiate or change studies, or stop the meta-analysis all-together, while still retaining Type-I error guarantees throughout. 

Without using the ALL-IN terminology, \cite{duan2020interactive} design several martingale methods to sequentially test a global null when each study ends with a p-value (instead of e-value) that is valid conditional on all past studies.  These new methods can be seen as sequential analogs to several well known nonsequential p-value combination rules like Fisher's or Stouffer's. Alternatively, one could \emph{calibrate} the p-values into e-values (calibrators are defined in~Section~\ref{sec:eci}) and multiply them as done above.

\subsection{False Discovery Rate}\label{subsec:fdr}

The false discovery rate (FDR) is probably the most popular error metric in modern large-scale multiple testing. The BH procedure~\citep{benjamini1995controlling} is the standard procedure for controlling the FDR when working with p-values. Given a target FDR level $\alpha$, it proclaims as discoveries the hypotheses corresponding to the $k^*$ smallest p-values out of $K$, where
\[
k^*:=\max\{k\in\{1,\dots,K\}: p_{(k)} \leq \alpha k / K\},
\] 
and $p_{(k)}$ represents the $k$-th smallest p-value. It is known to control the false discovery rate when the null p-values are independent of each other and of the non-nulls, as well as under a particular type of positive dependence known as PRDS~\citep{benjamini2001control}.

\cite{wang2022false} define an analogous e-BH procedure, which rejects hypotheses corresponding to the $k^*$ largest e-values, where
\[
k^*:=\max\{k\in\{1,\dots,K\}: e_{[k]} \geq  K/(k\alpha)\},
\] 
and $e_{[k]}$ represents the $k$-th largest e-value. Surprisingly, this procedure controls the FDR at $\alpha$ under arbitrary dependence between the e-values. The same result holds if one picks any set of $S$ e-values that are all larger than $K/(S\alpha)$; an analogous result does not hold for p-values.

\cite{xu2021unified} extended these results to \emph{bandit multiple testing}. There, the data to test the $K$ hypotheses is not available in advance, but must be collected adaptively, for example by assigning later subjects to more promising treatments as revealed by the results on earlier subjects. For each of the $K$ treatments, one can form an e-process to test the null hypothesis that the treatment effect is nonpositive. The $K$ e-processes have a complex dependence structure because of the adaptive assignment mechanism. Nevertheless, at any data-dependent stopping time, the e-BH procedure applied to the stopped e-processes  controls the FDR.

When both p-values and e-values are available for the same set of hypotheses (for example, from different datasets collected under different conditions), \cite{wang2022values} define generalizations of the above procedures that use both sources of information. In particular, e-values can serve as \emph{unnormalized weights} within standard FDR methods that use weighted p-values. The waiving of the need to normalize the weights (to sum to one) gives the e-value weighted methods a distinct power advantage over the usual normalized weights that are employed in weighted multiple testing.

\subsection{False Coverage Rate}
\label{sec:eci}
Suppose data regarding $K$ parameters has been collected, a data-dependent selection rule $\mathcal S$ is applied to select a subset $S$ of the parameters deemed of interest, and CIs for the selected parameters must be reported so as to keep the expected fraction of miscovering intervals at $\alpha$. The BY procedure \citep{benjamini2005false} is an analog of the BH procedure for this task: we report $(1-\alpha R/K)$-CIs for the selected parameters, where $1 \leq R \leq \lvert S\rvert$ is some function of the selection rule and dependence structure. Under certain dependence assumptions, this is proven to control the FCR at level $\alpha$. 

In contrast, the e-BY procedure of~\cite{xu2022post} applies only to e-CIs, which are CIs constructed by inverting tests based on e-values (discussed next). The authors prove that reporting $(1-\alpha \lvert S\rvert/K)$ e-CIs controls the FCR at $\alpha$ for any dependence structure, and any data-dependent selection rule $\mathcal S$ (including one that is fully aware of the corrected intervals).

\if
Formally, using the notation of Section~\ref{sec:confidence}, e-confidence sets for data $X^{\tau}$ can be defined whenever we are given a family of e-variables ${\cal E}:= \{ \prode{\tau}{\theta}: \theta \in \Theta \}$, where $\theta$ is a parameter of interest, and for all $\theta \in \Theta$, $\prode{\tau}{\theta}$ is an e-variable defined relative to null $\cP_{\theta}= \{P \in \Pi: \phi(\Pi) = \theta\}$.
Then the $(1-\alpha)$-e-confidence set based on family ${\cal E}$ is simply given by $\{ \theta \in \Theta: \prode{\tau}{\theta} < 1/\alpha\}$; by Markov's inequality, this is also a confidence set in the usual sense. If the e-confidence set is an interval, we refer to it as an e-CI. 
\fi

Concrete examples of e-CIs include all confidence sets based on universal inference, any (arbitrarily) stopped confidence sequence, and CIs constructed using Chernoff-style concentration inequalities. Further,~\cite{xu2022post} show that any CI can be converted to an e-CI by calibration. A calibrator $f$ is nonincreasing function from $[0,1]$ to $[0,\infty)$ such that $\int_0^1 f(x)dx = 1$. If a calibrator $f$ is also continuous at $1/\alpha$, then any CI constructed at (the more stringent) level $\alpha' := f^{-1}(1/\alpha)$ is an e-CI at level $\alpha$. This e-CI is always larger than the original CI.

As before, the implications for bandit multiple testing are interesting. One can construct and continuously monitor a CS for the effect size of each treatment, decide when to stop adaptively, select any subset for further study, and report corrected CIs using the stopped CSs at the e-BY adjusted level. As one example, one could run e-BH continually using the underlying e-processes, decide when to stop based on its rejections; then the corrected CIs will be congruent with the reported discoveries in the sense that all the corrected CIs will not contain the null parameter and both FDR and FCR will be controlled at level $\alpha$. 

\if
(The congruence is a result of the duality between tests and CIs that we employed earlier: the e-process exceeds $1/\alpha$ if and only if the CS does not intersect the null.)
\fi

\subsection{The Inevitability of e-hacking}

Peeking at the data obtained so far in order to decide whether to continue is only one of many abuses of statistical testing that have been classified as ``p-hacking''.  Because of the interpretation in terms of betting, peeking is legitimate when we test by e-values, but other abuses are neither legitimized nor prevented.  When statisticians commit these abuses using e-values, they have merely replaced ``p-hacking'' with ``e-hacking''.  A statistician is e-hacking, for example, whenever they implement many betting strategies with given data and report only the one that yields the greatest wealth.

The fundamental principle of testing by betting is that a bet on an outcome must be made before the outcome is observed.  Optional continuation is allowed in the case of successive bets because this condition is still met for each individual bet.  But claiming you would have bet in a certain way after you know the outcome is still humbug.  We can only hope that the clarity of these principles, even for laypeople, will make the possibilities for abuse more obvious and increase the pressure to distinguish between exploratory and confirmatory analysis.

In some situations, abuses can be prevented or mitigated by a separation of roles.  The use of e-values rather than p-values may be helpful in these situations.  

In academic disciplines where abuses may be driven by the need to publish, for example, editors can encourage preregistration of a study's data collection and analysis.  If we agree that the analysis should use the betting strategy that maximizes the expected logarithm of wealth under a reasonable alternative, then the proposed analysis necessarily identifies the alleged reasonable alternative; \cite{Shafer:2021} calls this the \emph{implied alternative}.  Editors and referees could reject proposed registrations for which this implied alternative is not really plausible and even agree in advance to publish the study when it is plausible and interesting.  This option does not arise when classical significance testing is used, because usually there is no unique alternative for which a test is most powerful.

When the statistician is embedded in a larger scientific enterprise, decisions about each step in data collection can be the result of consultation between the statistician and other scientists.  In the first flush of excitement about Wald's sequential analysis, \cite{Barnard:1947} saw this as the future of statistics, but it has been in tension with the notion of a p-value based on a global test statistic.  Testing by betting escapes this tension and  can be used even in collaborative meta-analysis (Section~\ref{sec:meta}).

\section{Other applications}

Game-theoretic statistics is rapidly evolving. Here are additional topics where it is relevant.

\subsubsection*{Comparing/Evaluating Forecasters.} Many experts and pundits now repeatedly make predictions about the weather, wars, sport games, business events, and elections probabilistically, sometimes as the probability of an event (one team beating another) or a predictive distribution (over the amount of rain the next day). How can we test whether probabilistic forecasters are doing a good job (are calibrated, for example), and how can we compare two different probabilistic forecasters? Such questions have been addressed in a game-theoretic setup by several recent works that use test supermartingales~\citep{henzi2021valid,arnold2021sequentially} or e-processes and confidence sequences~\citep{choe2021comparing}.

A fascinating general phenomenon, called \emph{Jeffreys's law} by \citet[Sect.~5.2]{Dawid84} in honor of Harold Jeffreys, is that two reliable forecasters must agree in the long run: if they differ too much, a Skeptic observing both of them will be able to discredit at least one of them \cite[Sect.~10.7]{Shafer/Vovk:2019}.

\subsubsection*{Multi-Armed Bandits and Reinforcement Learning.}
In sequential decision making, as modeled by a contextual multi-armed bandit or a reinforcement learning problem, one sees a sequence of ``contexts'' $x_t \in \mathcal X$, and one must decide which action $a_t \in \mathcal A$ to take in order to maximize a (discounted) sum of observed rewards $R(x_t,a_t)$. A policy $\pi$ is a mapping from $\mathcal X$ to $\mathcal A$, and one usually attempts to understand the unknown reward function $R$ by playing some exploratory policy $\pi_0$. One central question is the following: if the data was collected using some $\pi_0$, is it possible to estimate the quality (called ``value'') of some other policy $\pi_1$ that was never deployed? This is called ``off-policy evaluation'', and is a central problem of great practical interest. Recently, \cite{karampatziakis2021off} developed confidence sequences for off-policy evaluation when the rewards are bounded, and extended by~\cite{waudby2022anytime} to settings with unbounded importance weights and time-varying policies. We remark that outside of the off-policy setting, CSs are commonly used in contextual bandits~\citep{abbasi2011improved,chowdhury2017kernelized} and best arm identification~\citep{jamieson2014lil,kaufmann2021mixture}.

\section{Discussion}

\subsection{Connections with Other Areas}
\label{sec:relationship}

\subsubsection{Bayesian and Evidentialist Approaches.}
\label{subsubsec:bayes}
We already alluded to various connections between Bayesian and game-theoretic statistics (Bayes factors, Section~\ref{subsubsec:bfactor}, Jeffreys' prior, Section~\ref{sec:compositeH1}, Bayesian t-test Section~\ref{sec:ttest}). Even though interpretations are very different, a precise comparison would fill up an entire paper. 
We do highlight some relations in Appendix~\ref{app:bayes},  emphasizing that e-processes generalize likelihood ratios and, like these, can be interpreted as {\em evidence}. In the present section, we restrict ourselves to one specific simple technique to construct CSs, the {\em prior-posterior ratio martingale\/} ~\citep{waudby2020confidence,Grunwald23} so as to demonstrate  how Bayesian tools can be transformed into SAVI tools.
Suppose the data are drawn from $P_{\theta^*}$ for some unknown $\theta^*\in\Theta$ (extension to Bayesian nonparametrics is straightforward~\citep{neiswanger2021uncertainty}). Let $\pi_0 (\cdot)$ be a ``prior'' distribution over $\Theta$; we call this the working prior, because no assumptions are made about it.
After seeing $X_1,X_2,\dots,X_t$, let $\pi_t(\cdot)$ be the posterior distribution obtained via Bayes' rule. The central observation is that the density ratio $d\pi_0(\theta^*)/d\pi_t(\theta^*)$ is a test martingale for $P_{\theta^*}$, termed the ``prior-posterior ratio martingale'' (it is used for different purposes in Bayesian statistics, where it is called the {\em Dickey-Savage ratio}.  %
Thus, $\{\theta \in \Theta: d\pi_0(\theta)/d\pi_t(\theta) < 1/\alpha\}$ is a $(1-\alpha)$-CS for $\theta^*$.
Intriguingly, quite recently it has also been suggested within the Bayesian community  \citep{wagenmakers2020support,PawelLW22}  to use this CS, when applied to fixed $t$, as an alternative for the Bayesian posterior credible interval;  we comment further in Appendix~\ref{app:bayes}. These are also closely related to Bayesian ``snug regions'' proposed\footnote{We thank an anonymous reviewer for this reference.} by~\cite{hildreth1963bayesian}.

\subsubsection{Group-Sequential and Alpha Spending Methods.}\label{sec:groupseqs} 
These methods are mostly used in the clinical trial literature. Like our methods, they have their roots in the work of Robbins, Siegmund, Lai and others on anytime-valid tests in the 1970s. But they developed in quite a different direction: although there are exceptions\footnote{We thank J. Goeman and J. ter Schure for pointing this out to us.}  such as \cite{HuCG07}---group sequential methods provide Type-I error control under multiple looks at the data, but they typically require a prespecified final sample size, and a prespecified set of times at which one looks at the data. In contrast, e-processes can be updated as long as new data is available; an extensive comparison in the setting of the logrank test is performed by \cite{Schure/Grunwald/Ly:2021}. Principles for designing e-processes, such as the GRO criteria, or the RIPr and UI methods, do not seem to have analogues in the $\alpha$-spending/group sequential literature. But a firmer understanding of connections is desirable.

\subsubsection{Information Theory and Online Learning.} 
\label{subsubsec:mdl}
We touched on the relationship between our methods and the information-theoretic \emph{Minimum Description Length (MDL)} paradigm for model selection, learning and prediction \citep{BarronRY98,GrunwaldR20}  when discussing the REGROW criterion in Section~\ref{sec:compositeH1}. The connection to MDL and the related idea of universal coding runs quite deeply, due to Kraft's inequality, which states that for any probability distribution $\bar{Q}$ with probability mass function $\bar{q}$ and any stopping time $\tau$, there is a lossless code such for every realization $x^{\tau}$, the codelength achieved with this code is equal, up to a negligible roundoff term, to $- \log \bar{q}(x^{\tau})$; conversely, for any lossless code there is a distribution $\bar{Q}$ such that this correspondence holds. In MDL approaches one proceeds by associating statistical models (sets of distributions) ${\cQ}$ with `universal codes', represented as distributions $\bar{q}$ such that the codelengths are $-\log \bar{q}(x^t)$, designed to give small codelengths to the data at hand whenever the code corresponding to any element $P\in {\cQ}$  assigns a small codelength to the data. This very closely mirrors the construction of $\bar{q}$ via an estimator $\hat\theta$ or the method of mixtures as in Section~\ref{sec:compositeH1}. MDL model selection between a number of parametric models ${\cQ}_{\gamma}, \gamma \in \Gamma$ works by first associating each ${\cQ}_{\gamma}$ with a $\bar{q}_{\gamma}$ as above, and then picking as `the best explanation for data $x^t$' the $\gamma$ for which the associated codelength $- \log \bar{q}_{\gamma}(x^t)$ is minimal, reporting as evidence of model ${\cQ}_{\gamma_1}$ over ${\cQ}_{\gamma_2}$ the codelength difference $- \log \bar{q}_{\gamma_2}(x^t) - [- \log \bar{q}_{\gamma_1}(x^t)]$.  As a result, if there are just two models, $\gamma \in \{0,1\}$ and the null model is simple, the MDL approach is essentially equivalent to doing a test between null ${\cQ}_0$ and alternative ${\cQ}_1$ and reporting as evidence the logarithm of the e-value $\bar{q}_{1}(x^t)/q_0(x^t)$ \citep{GrunwaldR20} --- exactly the same as in Section~\ref{sec:compositeH1} but with evidence expressed on a logarithmic scale. When the null is composite, MDL and SAVI methods diverge, but we conjecture that  e-processes have a codelength interpretation---but with different codes than in classical MDL approaches.

One may also think of `universal codes' $\bar{q}$ as sequential prediction strategies that predict $x_t$ using $q(x_t \mid x^{t-1})$ and with loss assessed by the {\em logarithmic loss function\/} $-\log q(x_t \mid x^{t-1})$. The vast field of {\em online learning\/} is about such sequential prediction and the logarithmic loss takes an important special place in it.  Not surprisingly then, sequential prediction strategies from the online learning literature can often be  converted to  provide good (in some cases optimal in some sense) betting strategies for several problems. These connections have been emphasized by~\cite{orabona2021tight,waudby2020estimating,shekhar2021game,ramdas2022testing,casgrain2022anytime} and others.
Importantly, whereas in this paper we emphasized the case that $\cQ$ is a class of alternatives that are seriously contemplated as having generated the data, we may also do our tests with $\bar{q}$ that we suspect will predict the data reasonably well but cannot be considered `potentially true' in any meaningful sense. For example, the $\bar{q}$ may just be an {\em expert}, a probabilistic predictor, the inner workings of which may be completely unknown to us. 

\subsection{Open Questions}
\label{sec:open}


\subsubsection{Existence of e-processes.}  For what classes of distributions $\cP$ are there nontrivial (a)  test martingales, (b) test supermartingales but no test martingales, (c) e-processes but no test supermartingales, (d) none of the above? While \cite{ramdas2022testing,ruf2022composite} have interesting examples separating the above concepts, this separation is not yet fully understood in general. A recent preprint by \cite{zhang2023exact} yields new insights via convex geometry and optimal transport.

\subsubsection{Choice of Filtration.}
\label{subsubsec:choice}

The choice of a filtration is a design choice, and one need not choose the richest one, the one generated by the observations.  The choice 
affects both
safety (the set of stopping times $\tau$ for which the expected value of $E_\tau$ does not exceed one under the null) and power (how fast the wealth grows under the alternative).

Recall the example of testing exchangeability in Section~\ref{sec:exch+LC}. As explained, there is no nontrivial test martingale for the problem in the filtration of the observations, but there is one in the coarsened filtration of conformal p-values. Despite the generality of conformal p-values, the coarsening implies safety under optional stopping for a smaller set of stopping times that cannot see the original data, but only the conformal p-values.  This sacrifice is unnecessary for data with a known discrete support, in which case an e-process is available in the original filtration, that appears in experiments to be at least as powerful as the conformal test martingale.

The picture is a little different for the test martingales in a shrunk filtration constructed by \cite{perez2022estatistics} for the problem of testing group-invariant hypotheses (such as the t-test example of \cite{Cox52} and \cite{Lai:1976} discussed in Section \ref{sec:ttest}).  These are not e-processes with respect to the original filtration, but are in the shrunk one, and thus have a weaker guarantee under the null,
 but they still maximize the rate of growth amongst all e-processes, even those with respect to the original filtration. Thus, they have worse safety properties but better growth than (say) universal inference. 

When can one coarsen the filtration in order to design useful new e-processes, and when are these more or less powerful than ones in the original filtration?

\subsubsection{Admissibility.}
\label{subsec:admissibility} Can one characterize admissibility of an e-process succinctly, with a condition that is both  necessary and sufficient? \cite{ramdas2020admissible} provide both necessary and sufficient conditions for admissibility, but currently these do not match. For example, they prove that, if there exists a common dominating measure, then $E$ being admissible implies that $E_t = \inf_{P \in \cP} M^P_t$, where $M^P_t$ is a test martingale for $P$. The universal inference e-process has this form. But this condition is not sufficient for admissibility: e-processes satisfying $E_t = \inf_{P \in \cP} M^P_t$ may not always be admissible (indeed, universal inference has this form, and we know examples where it is inadmissible). For admissibility, the $\{M^P_t\}_{P \in \cP}$ need to agree to some extent---they need to be large or small on similar events; if on each event, some test martingales are large while others are small, the infimum will always be small. Given that admissibility is a low bar, delineating this need for agreement is an important open problem. Intriguingly, the use independent external randomization can significantly alter the story~\citep{ramdas2023randomized}.

\subsubsection{Questions about RIPr.} When exactly does the RIPr procedure applied to data $X^t$ separately for each sample size $t$ yield an e-process, as opposed to just a sequence of e-variables? (Section~\ref{sec:ripr}). When the RIPr yields an e-process, there is strong justification to use it, but how much is lost if it is replaced by the (always applicable) universal inference e-process?  Understanding the power of universal inference is itself quite open; progress was made in the Gaussian setting by~\cite{dunn2021gaussian}.
In current applications the GRO-optimal RIPr e-variables are sometimes given by simple, analytic formulas (Section~\ref{sec:independence} and~\ref{sec:ttest}), but for other applications numerical optimization is required (e.g.\ the logrank test as in Section~\ref{sec:monotone}). An algorithm by \cite{Li99} is slow. Do there exist practically effective algorithms? 

\smallskip

We end by pointing the reader to Appendix~\ref{app:price}, where we discuss the question of whether there is a ``price'' to be paid by SAVI methods. We anticipate more discussion around this (part philosophical, part technical) topic.

\addcontentsline{toc}{section}{References}

\bibliographystyle{imsart-nameyear} 
\bibliography{ReferencesAPVG}

\appendix

\section{SAVI as a Frequentist -- Evidential -- Bayesian middle ground?}\label{app:bayes}


E-processes can be seen as quantifying evidence against a null hypothesis and are quite meaningful even without being used for a sequential test required to have some error probability, and even in a batch setting such as the multiple testing settings above. E-processes lack some of the properties of p-values that make the latter  less suitable to think of as `evidence' (such as the p-value's  dependency on whether or not particular actions are taken in counterfactual situations, as exemplified by examples in the literature such as Pratt's volt-meter story \citep{edwards1992likelihood}) and generalize the likelihood ratio that is embraced by the likelihoodists as the `right' formalization of relative evidence \citep{royall1997statistical}. 

Comparing our methods to Bayesian ones, we see that, with simple nulls, admissible e-processes and Bayes factors always coincide; in parametric tests with composite nulls, e-processes and Bayes factors sometimes (e.g.\ in the group invariant setting of Section~\ref{sec:ttest}) but not always coincide; and with nonparametric tests they start differing quite a lot. If it comes to confidence sequences and e-confidence intervals, we find that even in one-dimensonal parametric settings, $(1-\alpha)$- e-confidence intervals (and equivalently stopped confidence sequences) do not coincide with Bayesian $(1-\alpha)$-posterior credible intervals, the latter being significantly narrower. To see this, note that, from Section~\ref{sec:compositeH1} (and defining e-CIs as in Section~\ref{sec:eci}), we find that, for any fixed prior density $w$ on $\Theta \subset {\mathbb R}$, the family of e-variables $\{\prode{\tau}{\theta}: \theta \in \Theta \}$ for data $X^{\tau}$ with 
\[
\bar{P}(\theta \mid X^\tau) := \frac{1}{\prode{\tau}{\theta}} = 
\frac{p_{\theta}(X^{\tau})}{\int p_{\theta'}(X^{\tau}) w(\theta') d \theta'}
\]
(this is just the reciprocal of the  prior-posterior ratio martingale of Section~\ref{subsubsec:bayes} stopped at time $\tau$)
defines an e-confidence interval at level $(1-\alpha)$ as $\{\theta: \bar{P}(\theta \mid \tau) > \alpha \}$, whenever the latter set is an interval. 
By Bayes' theorem, the Bayes posterior based on the same prior $w$ is given by
\begin{equation*}
w(\theta \mid X^{\tau}) =  \frac{w(\theta) \cdot p_{\theta}(X^{\tau})}{\int p_{\theta'}(X^{\tau}) w(\theta') d \theta'} = 
w(\theta) \cdot \bar{P}(\theta \mid X^{\tau}),
\end{equation*}
and defines a posterior credible interval at level $(1-\alpha)$ as $[\theta_{L},\theta_R]$ chosen so that 
$$\E_{\theta \sim W} [{\bf 1}_{\theta \in [\theta_L,\theta_R]} \cdot \bar{P}(\theta \mid X^\tau)] = \int_{\theta_L}^{\theta_R} w(\theta \mid X^{\tau}) = 1-\alpha.$$ 
We see that all elements $\theta$ of an  e-confidence interval must have $\bar{P}(\theta \mid X^{\tau}) \geq \alpha$; for a Bayesian credible interval this only has to hold in average over the prior, causing the latter to be narrower in practice. 

Intriguingly though, some Bayesian statisticians have noted that the standard Bayesian posterior credible interval has no clear `evidential' interpretation. They instead propose a {\em Bayesian support interval\/} \citep{wagenmakers2020support}---also called `evidential support interval'---  where the $k$- support interval is the interval containing all parameter values under which the observed data $X^{\tau}$ are at least $k$ times as likely than under the Bayesian marginal distribution'. As \cite{PawelLW22} note, for simple nulls and $k < 1$, this actually coincides precisely with the $(1-k$)-e-confidence interval based on the family of e-values based on the same prior density $w$ as the Bayes marginal. 

We would venture that, for models $\Pi$ with  parameter of interest $\theta = \phi(P)$ and additional nuisance parameters, and also in nonparametric settings, the e-confidence intervals based on a family $\{\prode{\tau}{\theta} : \theta \in \Theta \}$ still have an `evidential' interpretation, although in this case they will usually not be equal to a Bayesian support interval any more. 

Taking the `e-values are similar to, but different from Bayes factors' line of reasoning even further, one could daringly suggest to define $\bar{P}(\theta \mid X^{\tau}) := 1/\prode{\tau}{\theta}$ as an analogue of the Bayesian posterior or confidence distributions, even for multiparameter and nonparametric problems in which it does not coincide with the Savage-Dickey density ratio.  
This was done informally in \citet[Appendix E7]{waudby2020confidence} who visualize uncertainty by drawing $\bar{P}(\theta \mid X^{\tau})$ as a function of $\theta$.
\cite{Grunwald22,Grunwald23} shows that this {\em e-posterior\/} can be motivated not just evidentially, but also decision-theoretically. Just like the Bayes posterior can, assuming the prior was chosen well, be used to obtain optimal decisions for arbitrary loss functions by combining posterior and loss in a certain way (minimizing Bayes-posterior expected loss), the e-posterior can be used as a basis for obtaining decisions with minimax optimality guarantees for arbitrary loss functions. The guarantees hold irrespective of the chosen prior, but become weaker the more atypical the data look with respect to the prior. 
In the same papers (see also \cite{BatesJSS22}), it is shown that, even in a nonsequential context,  standard Neyman-Pearson testing is not adequate  if the decision problem at hand (e.g.\ choose between the four actions $\{$ vaccinate no-one; only adults; only the elderly; or everyone$\}$) has more actions than just the `reject' and `accept' of the Neyman-Pearson theory; with a decision rule based on e-variables one can effectively deal with such---realistic---settings. This has direct repercussions for the reproducibility crisis: the standard Neyman-Pearson based approaches may simply not be suitable for the complex real-world problems that we apply our test results to. 

In conclusion, let us stress that we do not view the above observations as disqualifying the Bayesian, evidential or Neyman-Pearsonian paradigm. Rather, we feel that SAVI methods effectively unify some of the fundamental ideas of each; respectively: they allow to infuse prior knowledge into one's procedures;  they output numbers with a clear evidential meaning; and they ensure error control and coverage.

\section{Does SAVI come at a price?}
\label{app:price}

We often receive questions like: ``does the much greater flexibility allowed by SAVI methods compared to traditional Neyman-Pearson testing come at a price? Is there, for example, a loss of power or a need for larger samples before conclusions can be drawn? How competitive are SAVI methods with classical ones?''

Let us not be too starry-eyed here. There is always a fundamental information-theoretic tradeoff between different types of errors that is unavoidable by any framework, including ours. 

At one level, one may view game-theoretic statistics as trading off the two errors differently from classical statistics --- we ask for a more stringent form of ``type-I error'' that even holds under optional stopping or continuation of experiments, but accept a (slightly) higher ``type-II error'' in return, and many of the cited papers do detailed analysis on the quantitative nature of the aforementioned qualitative tradeoff, e.g.\ \cite{GrunwaldHK19,Schure/Grunwald/Ly:2021,waudby2020estimating,howard2021time,wang2022false}. 

However, the errors are in quotes because even our definitions of type-I and type-II errors themselves differ from classical statistics --- instead of minimizing the probability of each error, when constructing e-statistics (e-values,  test supermartingales, e-processes, etc.), the ``type-I'' error control corresponds to making sure the expectation of the e-statistic is at most one under the null, and low ``type II'' error corresponds to finding an e-statistic whose expected logarithmic value under the alternative is as large as possible.

Even talking about the ``price'' assumes that one set of metrics dominate the other, but we contend that our performance metrics are in fact more suitable in many situations,  in that the resulting framework has a large number of benefits for reproducibility --- even above the optional stopping/continuation, the methods are very robust to dependence (Section~\ref{subsec:fdr}), it is easy to combine evidence from independent or dependent studies (Sections~\ref{subsec:avg-e},~\ref{subsec:mult-e} and~\ref{sec:meta}), etc.
Nevertheless, many SAVI papers directly and extensively compare `real' Type-I and Type-II errors obtained with our approach (although that is not what we optimize for with the SAVI approach) to those obtained by classical frequentist approaches.

As to a comparison to Bayesian approaches, it is also difficult and of limited consequence to directly compare numbers, but in Appendix~\ref{app:bayes} we do give some feel about how Bayesian credible intervals compare to our anytime-valid confidence sets/intervals (which are wider than Bayesian credible intervals, but not wider than the support interval).

It is true that, if we compare an e-variable with GRO status for testing a specific $H_0$ and $H_1$ at a fixed sample size to the uniformly most powerful Neyman-Pearson test at that sample size (assuming it exists), the latter has more power. The exact difference (as well as, relatedly, difference of width between standard confidence intervals and anytime-valid confidence sequences, which tend to get between 1.5 and 2 times as wide) is investigated in detail, both theoretically and by simulation by the earlier cited papers (and many others).

   On the other hand, if one allows for optional stopping and stops as soon as one can reject the null, then {\em expected\/} minimal stopping times under the alternative are about the same or even smaller than the fixed $n$ needed to get the same power with a Neyman-Pearson test --- illustrating that it all depends how one  measures performance quality. 
   
   This suggests that (and indeed we feel that) really, the wrong questions are often being asked: the SAVI methods are optimizing for different criteria --- a measure of evidence (wealth) that grows fast under any alternative, and does not grow under any null 
   --- which arguably are often more suitable in the applied sciences in which replicability issues abound. Thus, comparison in terms of power is only of limited use. 

   \cite{Grunwald22,Grunwald23} also shows that based on e-processes, one has the flexibility to make decisions based on arbitrary loss functions that are determined post-hoc, which is arguably highly relevant for practice yet also is not captured by power.


\end{document}